\documentclass{elsarticle}

\DeclareMathAlphabet{\mathpzc}{OT1}{pzc}{m}{it}

\usepackage{bm}
\usepackage{graphicx}
\usepackage{amssymb}
\usepackage{amsmath}

\usepackage{natbib}
\usepackage{amsthm}
\usepackage[final]{showkeys}
\usepackage[hmargin=2.5cm,vmargin={3.5cm,4cm}]{geometry}

\numberwithin{equation}{section}
\allowdisplaybreaks

\newtheorem{thm}{Theorem}[section]
\newtheorem{rem}{Remark}[section]

\title{Randomly Stopped Nonlinear Fractional Birth Processes\tnoteref{t1}}

\author[eo]{Enzo Orsingher\corref{cor1}}

\author[fp]{Federico Polito}

\cortext[cor1]{Corresponding author.}
\fntext[fn2]{Dipartimento di Scienze Statistiche, ``Sapienza'' Universit\`a di Roma, Piazzale Aldo Moro 5, 00185 Rome, Italy. Piazzale Aldo Moro 5, 00185. Tel: +39 06 49910585, fax: +39 06 4959241. E-mail address: enzo.orsingher@uniroma1.it.}

\fntext[fn3]{Dipartimento di Scienze Statistiche, ``Sapienza'' Universit\`a di Roma, Piazzale Aldo Moro 5, 00185 Rome, Italy. Piazzale Aldo Moro 5, 00185. Tel: +39 06 49910499, fax: +39 06 4959241. E-mail address: federico.polito@uniroma1.it.}

\tnotetext[t1]{2010 Mathematics Subject Classification. Primary 60G22, 60G55.}

\begin{document}

	\begin{abstract}
		We present and analyse the nonlinear classical pure birth process $\mathpzc{N} (t)$, $t>0$,
		and the fractional pure birth process $\mathpzc{N}^\nu (t)$, $t>0$, subordinated to
		various random times,
		namely the first-passage time $T_t$ of the standard Brownian motion $B(t)$, $t>0$,
		the $\alpha$-stable subordinator $\mathpzc{S}^\alpha(t)$, $\alpha \in (0,1)$, and others.
		For all of them we derive the state probability distribution $\hat{p}_k (t)$,
		$k \geq 1$ and, in some cases, we also present the corresponding governing
		differential equation.

		We also highlight interesting interpretations for both the subordinated classical birth process
		$\hat{\mathpzc{N}} (t)$, $t>0$, and its fractional counterpart $\hat{\mathpzc{N}}^\nu (t)$, $t>0$
		in terms of classical birth processes with random rates evaluated
		on a stretched or squashed time scale.

		Various types of compositions of the fractional pure birth process $\mathpzc{N}^\nu(t)$
		have been examined in the last part of the paper. In particular, the processes
		$\mathpzc{N}^\nu(T_t)$, $\mathpzc{N}^\nu(\mathpzc{S}^\alpha(t))$, $\mathpzc{N}^\nu(T_{2\nu}(t))$,
		have been analysed, where
		$T_{2\nu}(t)$, $t>0$, is a process related to fractional diffusion equations.
		Also the related process $\mathpzc{N}(\mathpzc{S}^\alpha({T_{2\nu}(t)}))$ is investigated
		and compared with $\mathpzc{N}(T_{2\nu}(\mathpzc{S}^\alpha(t))) = \mathpzc{N}^\nu
		(\mathpzc{S}^\alpha(t))$. As a byproduct
		of our analysis, some formulae relating Mittag--Leffler functions are obtained.
	\end{abstract}

	
	\maketitle
	
	\section{Introduction}

		We here consider the pure birth process $\mathpzc{N}(t)$, $t>0$, (linear and nonlinear)
		composed with different processes like the first-passage time of Brownian motion $T_t$
		(possibly iterated n-times), the sojourn time of Brownian motion $\Gamma_t$ and bridge
		$\mathfrak{G}_t$, and $\alpha$-stable
		processes $\mathpzc{S}^\alpha(t)$.

		The subordination of processes (first introduced by \citet{bochner}) has been studied by
		several authors, over the years, in connection, for example, to modelling the wear of instruments
		during the real working time, or security trading which takes into account
		fluctuations of the economic activity during the time elapse $t$ (see \citet{lee}).

		The second part of the paper concerns the subordination of the fractional pure birth
		process $\mathpzc{N}^\nu(t)$, $t>0$, $0<\nu \leq 1$, with the processes
		$\mathpzc{S}^\alpha(t)$ and $T_{2\alpha}(t)$, establishes that
		$\mathpzc{N}^\nu(\mathpzc{S}^\alpha(t)) = \mathpzc{N}(T_{2\nu}(\mathpzc{S}^\alpha(t)))$,
		and discuss its connection with $\mathpzc{N}(\mathpzc{S}^\alpha(T_{2\nu}(t)))$.

		Subordinated processes connected with fractional and higher order partial differential equations
		are treated in numerous recent papers. Most of them concern compositions of time-continuous
		processes (see for example \citet{baeumer2}), but also point processes
		(\citet{laskin}, \citet{scalas}, \citet{cahoy2}, \citet{orsbeg2}, \citet{meer}).

		Birth processes stopped at different random times can be useful to model branching
		processes under laboratory conditions. For diseases started off artificially,
		the spread of the infected population can be stopped when the experiment leads to
		convincing conclusions. The cost of the investigation can play a certain role in stopping
		the artificially constructed experiment.
		The fluctuations of the temperature during the effective time $t$ can influence the growth
		rapidity of cells or of bacteria and thus the population size
		can be thought as a function of the temperature modelled as a random time process.
		The same reasoning underlies experiments
		in physical studies on chain reactions. In the case of $\mathpzc{N}(\Gamma_t)$,
		where $\Gamma_t$ is the sojourn time of a Brownian motion on the positive half-line,
		the experiment can be interrupted immediately (if it proves useless), or at the end
		of the time interval $[0,t]$ (in the case that no evidence can be attained in a short time).

		We recall that the distribution of the nonlinear fractional birth process (with one
		progenitor) reads
		\begin{equation}
			\label{nlinearnuintro}
			\text{Pr} \left\{ \mathpzc{N}^\nu(t) = k \mid \mathpzc{N}^\nu(0)=1 \right\} =
			\begin{cases}
				\prod_{j=1}^{k-1} \lambda_j
				\sum_{m=1}^k \frac{ E_{\nu,1} (-
				\lambda_m t^\nu)}{\prod_{
				l=1,l \neq m}^k \left( \lambda_l -
				\lambda_m \right) }, & k > 1, \\
				E_{\nu,1} (-\lambda_1 t^\nu), & k=1,
			\end{cases}
		\end{equation}
		where
		\begin{equation}
			\label{mittagintro}
			E_{\nu,\gamma} \left( x \right) = \sum_{h=0}^\infty
			\frac{x^h}{\Gamma \left( \nu h+\gamma \right)},
		\end{equation}
		is the Mittag--Leffler function and $\lambda_k$, $k \geq 1$, are the birth rates (see \citet{pol}).

		For $\lambda_k=\lambda \cdot k$ (fractional linear birth process), formula \eqref{nlinearnuintro}
		takes the simple form
		\begin{equation}
			\label{linearintro}
			\text{Pr} \left\{ N^\nu(t) =k \mid N^\nu(0)=1 \right\}
			= \sum_{m=1}^k \binom{k-1}{m-1} (-1)^{m-1}
			E_{\nu,1} (-\lambda mt^\nu), \qquad k\geq1, \: t>0.
		\end{equation}

		For $\nu=1$, we retrieve from \eqref{nlinearnuintro} and \eqref{linearintro}
		the classical distributions of nonlinear and linear pure birth process, by taking
		into account that $E_{1,1}(x) = e^x$.

		The simplest subordinator considered is the first-passage time
		\begin{equation}
			\label{stop}
			T_t = \inf \left\{ s \colon B(s)=t \right\},
		\end{equation}
		where $B$ is a standard Brownian motion, independent of the birth process
		considered.
		For us it is relevant that the probability density of \eqref{stop}
		\begin{equation}
			q(t,s)ds=\text{Pr} \left\{ T_t \in ds \right\},
		\end{equation}
		satisfies the following equation
		\begin{equation}
			\label{taehintro}
			\frac{\partial^2}{\partial t^2} q(t,s) = 2 \frac{\partial}{\partial
			s} q(t,s), \qquad t>0, \: s>0,
		\end{equation}
		as a direct check shows.

		In view of \eqref{taehintro} we can establish the following relation between the state
		probabilities
		\begin{equation}
			\label{subintro}
			\hat{\mathpzc{p}}_k^\nu(t) = \text{Pr} \left\{ \mathpzc{N}^\nu (T_t)=k \right\}
		\end{equation}
		and \eqref{nlinearnuintro}:
		\begin{equation}
			\label{diffintro}
			\frac{d^2}{dt^2} \hat{\mathpzc{p}}_k^\nu(t) = -2 \int_0^\infty q(t,s)
			\frac{d}{ds} \text{Pr} \left\{\mathpzc{N}^\nu (s) =k \right\} ds.
		\end{equation}
		For $\nu=1$, equation \eqref{diffintro} becomes the second-order difference-differential
		equation
		\begin{equation}
			\label{eq-secintro}
			\frac{d^2}{d t^2} \hat{\mathpzc{p}}_k^\nu (t)
			= 2 \left[ \lambda_k \hat{\mathpzc{p}}_k^\nu (t) - \lambda_{k-1}
			\hat{\mathpzc{p}}_{k-1}^\nu (t) \right], \qquad k \geq 1.
		\end{equation}
		Furthermore, for $\nu=1$, the probability distribution \eqref{subintro} can be worked out
		explicitely and becomes
		\begin{equation}
			\label{nlinear-compintro}
			\hat{\mathpzc{p}}_k (t) =
			\begin{cases}
				\prod_{j=1}^{k-1} \lambda_j
				\sum_{m=1}^k \frac{ e^{-
				t \sqrt{2 \lambda_m}}}{\prod_{
				l=1,l \neq m}^k \left( \lambda_l -
				\lambda_m \right) }, & k > 1, \: t>0, \\
				e^{-t \sqrt{2 \lambda_1}}, & k=1, \: t>0.
			\end{cases}
		\end{equation}
		For $0<\nu<1$, in light of the well-known integral representation of the Mittag--Leffler
		function
		\begin{equation}
			\label{mitta-intintro}
			E_{\nu,1} (- \lambda t^\nu) = \frac{\sin \nu \pi}{\pi}
			\int_0^\infty \frac{r^{\nu-1} e^{-r \lambda^{\frac{1}{\nu}} t}}{
			r^{2 \nu} + 2 r^\nu \cos \nu \pi + 1} dr, \qquad \nu \in (0,1),
		\end{equation}
		we obtain several different representations of the distributions of the subordinated
		processes.

		For $\nu=1/2$, we have the following result
		\begin{equation}
			\label{interpintro}
			\hat{\mathpzc{p}}_k^{\frac{1}{2}} (t) =
			\frac{ \sqrt{2} }{ \pi } \int_0^\infty \frac{1}{ \left( \frac{w^2}{2}
			+1 \right) } \text{Pr} \left\{ \mathpzc{N}_w (t) = k \mid \mathpzc{N}_w
			(0) =1 \right\} dw,
		\end{equation}
		which shows that $\mathpzc{N}^{1/2}(T_t)$ is equivalent in distribution to a fractional
		pure birth process (denoted by $\mathpzc{N}_{W}(t)$)
		with rates $\lambda_k \cdot W$, where $W$ is a folded Cauchy
		distribution with scale parameter equal to $\sqrt{2}$.

		We have also that
		\begin{equation}
			\label{interp2intro}
			\hat{\mathpzc{p}}_k^{\frac{1}{2}} (t) =
			\int_0^\infty \text{Pr} \left\{ \mathpzc{N} (s) = k \right\}
			\text{Pr} \left\{ \bigl| C(\sqrt{2} t ) \bigr| \in ds \right\}.
		\end{equation}
		In other words, $\mathpzc{N}^{1/2} (T_t)$ is also equivalent in distribution
		to $\mathpzc{N}(| C(\sqrt{2}t) |)$, $C$ being a Cauchy process.

		We generalise the previous framework by considering the iterated process
		\begin{equation}
			\label{iterateintro}
			\tilde{\mathpzc{N}^\nu} (t) =
			\mathpzc{N}^\nu \biggl[ T^1_{T^2_{._{._{{T^n_t} }}}} \biggr],
			\qquad t>0,
		\end{equation}
		where $T^1_t, \dots, T^n_t$, are independent first-passage times and
		\begin{equation}
			T^j_{T^{j+1}_{._{._{{T^n_t} }}}} = \inf \left\{ s \colon B^j(s) =
			T^{j+1}_{T^{j+2}_{._{._{{T^n_t} }}}} \right\}, \qquad j=1,\dots,(n-1),
		\end{equation}
		where $B^j(t)$, $t>0$, $1\leq j \leq n$, are independent Brownian motions.
		In particular, for $\nu=1$ we show that the state
		probabilities
		\begin{equation}
			\text{Pr} \left\{ \tilde{\mathpzc{N}}^1(t)=k \right\} = \tilde{\mathpzc{p}}_k^1(t),
		\end{equation}
		satisfy the $2^n$th order equations
		\begin{equation}
			\frac{d^{2^n}}{dt^{2^n}} \tilde{\mathpzc{p}}_k^1(t) =
			2^{2^n-1} \left\{ \lambda_k \tilde{\mathpzc{p}}_k^1(t) -
			\lambda_{k-1} \tilde{\mathpzc{p}}_{k-1}^1(t)\right\}.
		\end{equation}
		The distribution $\tilde{\mathpzc{N}}^1(t) = \tilde{\mathpzc{N}}(t)$ (for short)
		is directly derived and reads
		\begin{equation}
			\label{nlinear-mcompintro}
			\tilde{\mathpzc{p}}_k (t) =
			\begin{cases}
				\prod_{j=1}^{k-1} \lambda_j
				\sum_{m=1}^k \frac{ e^{-t
				\lambda_m^{\frac{1}{
				2^n}} 2^{\left(1-\frac{1}{
				2^n} \right)}}}{\prod_{
				l=1,l \neq m}^k \left( \lambda_l -
				\lambda_m \right) }, & k > 1, \\
				e^{-t \lambda_1^{\frac{1}{
				2^n}} 2^{\left(1-\frac{1}{
				2^n} \right)}}, & k=1.
			\end{cases}
		\end{equation}
		For $n\rightarrow\infty$, we obtain from \eqref{nlinear-mcompintro} that
		\begin{equation}
			\lim_{n\rightarrow\infty} \tilde{\mathpzc{p}}_k (t)=
			\begin{cases}
				e^{-2t}, & k =1, \\
				0, & k > 1.
			\end{cases}
		\end{equation}

		In the last part of the paper we examine different types of compositions of the fractional
		pure birth process with positively skewed stable processes $\mathpzc{S}^\alpha (t)$, $t>0$,
		$0<\alpha\leq 1$. For $\alpha=\nu$, we show that
		\begin{equation}
			\mathpzc{N}^\nu (\mathpzc{S}^\nu (t))
			\overset{\text{i.d.}}{=}\mathpzc{N}(t\mathcal{W}_\nu), \qquad 0<\nu<1.
		\end{equation}
		For the stable random variables
		$\mathpzc{S}^\nu_1$, $\mathpzc{S}^\nu_2$, it is well-known that the ratio
		\begin{equation}
			\label{wallet}
			\mathcal{W}_\alpha =
			\left( \frac{ \mathpzc{S}^\nu_1 }{ \mathpzc{S}^\nu_2 } \right)^\alpha
		\end{equation}
		(sometimes called Lamperti law),
		has probability density equal to
		\begin{equation}
			f_{\mathcal{W}_\alpha}(r) = \frac{\sin \nu \pi}{\alpha \pi}
			\frac{r^{\frac{\nu}{\alpha}-1}}{r^{2\frac{\nu}{\alpha}} + 2r^\frac{\nu}{\alpha}
			\cos \nu \pi + 1}, \qquad r>0.
		\end{equation}

		Furthermore, we show that $\mathpzc{N}^\nu(\mathpzc{S}^\nu (t)) =
		\mathpzc{N}(T_{2\nu}(\mathpzc{S}^\nu(t))) \overset{\text{i.d.}}{\neq} \mathpzc{N}
		\left( \mathpzc{S}^\nu ({T_{2\nu}(t)}) \right)$. We are also able to prove that
		\begin{equation}
			\mathpzc{N}^\nu(T_{2\alpha}(t)) = \mathpzc{N} (T_{2\nu}(T_{2\alpha}(t)))
			\overset{\text{i.d.}}{=} \mathpzc{N} (T_{2\nu\alpha}(t)) =
			\mathpzc{N}^{\nu\alpha} (t).
		\end{equation}
		As a byproduct of our analysis we obtain the following integral relation between
		Mittag--Leffler functions of different indices:
		\begin{align}
			E_{\nu\alpha,1}(-\lambda_mt^{\nu\alpha}) & = \frac{\sin\nu\pi}{\pi}
			\int_0^\infty \frac{r^{\nu-1}}{r^{2\nu} +2r^\nu \cos \nu\pi+1} E_{\alpha,1}
			\left( -r \lambda_m^{\frac{1}{\nu}}t^\alpha \right) dr \\
			& = \frac{\sin\alpha\pi}{\pi}
			\int_0^\infty \frac{r^{\alpha-1}}{r^{2\alpha} +2r^\alpha \cos \alpha\pi+1} E_{\nu,1}
			\left( -r \lambda_m^{\frac{1}{\alpha}}t^\nu \right) dr, \qquad 0< \alpha,\nu \leq 1.
			\notag
		\end{align}

	\section{Subordinated nonlinear birth processes}

		In this section we study in detail the nonlinear pure birth process stopped at $T_t$ and
		we derive the state probabilities $\hat{\mathpzc{p}} = \text{Pr} \left\{ \mathpzc{N}
		(T_t) = k \mid \mathpzc{N}(0) = 1 \right\}$, $k\geq 1$, and the corresponding
		governing differential equations.

		We give some information about the process $\mathpzc{N}(t)$, $t>0$, evaluate explicitly
		its mean value $\mathbb{E}\mathpzc{N}(t)$, and discuss also the linear birth process
		(sometimes referred as Yule--Furry process).

		\subsection{Preliminaries}

			The state probabilities
			$\mathpzc{p}_k (t) = \text{Pr} \left\{ \mathpzc{N} \left(
			t \right) = k \mid \mathpzc{N}(0)=1 \right\}$
			read (see e.g.\ \citet{skoro}, page 322)
			\begin{equation}
				\label{nlinear}
				\mathpzc{p}_k (t) =
				\begin{cases}
					\prod_{j=1}^{k-1} \lambda_j
					\sum_{m=1}^k \frac{ e^{-
					\lambda_m t}}{\prod_{l=1, l \neq m}^k \left( \lambda_l -
					\lambda_m \right) }, & k > 1, \: t>0, \\
					e^{-\lambda_1 t}, & k=1, \: t>0.
				\end{cases}
			\end{equation}
			For the case of $n_0$ progenitors (see \citet{chiang}, page 51),
			formula \eqref{nlinear} must be replaced by
			\begin{equation}
				\label{nlinearn0}
				\mathpzc{p}_k (t) =
				\begin{cases}
					\prod_{j=n_0}^{k-1} \lambda_j
					\sum_{m=n_0}^k \frac{ e^{-
					\lambda_m t}}{\prod_{l=n_0, l \neq m}^k \left( \lambda_l -
					\lambda_m \right) }, & k > n_0, \: t>0, \\
					e^{-\lambda_{n_0} t}, & k=n_0, \: t>0.
				\end{cases}
			\end{equation}
			We assume that $\sum_k 1/\lambda_k = \infty$ in such a way that the process is
			non-exploding (see \citet{feller1}, page 452). For a discussion on this point,
			consult \citet{grimmett}, page 252.
			The probabilities \eqref{nlinear} satisfy the following difference-differential
			equations:
			\begin{equation}
				\label{fra}
				\frac{d}{dt} \mathpzc{p}_k (t) =
				- \lambda_k \mathpzc{p}_k(t) + \lambda_{k-1} \mathpzc{p}_{k-1}(t), \qquad k \geq 1.
			\end{equation}

			We have our first result in the next theorem.
			\begin{thm}
				\label{ball}
				The mean value of the nonlinear birth process is
				\begin{align}
					\label{valatt2}
					\mathbb{E} \mathpzc{N}(t) = 1 + \sum_{k=1}^\infty
					\left\{ 1 - \sum_{m=1}^k \prod_{l=1,l\neq m
					}^k \frac{\lambda_l}{\lambda_l-\lambda_m}
					e^{-\lambda_m t} \right\}.
				\end{align}
				\begin{proof}
					From equation \eqref{fra}, we have that
					\begin{align}
						\sum_{k=1}^\infty k \frac{d}{dt} \mathpzc{p}_k(t) & =
						- \sum_{k=1}^\infty k \lambda_k \mathpzc{p}_k(t) +
						\sum_{k=2}^\infty k \lambda_{k-1} \mathpzc{p}_{k-1}
						= \sum_{k=1}^\infty \lambda_k \mathpzc{p}_k(t).
					\end{align}
					By integrating both members in $(0,t)$, we obtain
					\begin{align}
						\label{levis}
						\sum_{k=1}^\infty k \mathpzc{p}_k(t) -1 & =
						\sum_{k=1}^\infty \lambda_k \int_0^t \mathpzc{p}_k(s) ds \\
						& = \lambda_1 \int_0^t \mathpzc{p}_1 (s) ds +
						\sum_{k=2}^\infty \lambda_k \left\{ \int_0^t
						\prod_{j=1}^{k-1} \lambda_j
						\sum_{m=1}^k \frac{ e^{-
						\lambda_m s}}{\prod_{l=1, l \neq m}^k \left( \lambda_l -
						\lambda_m \right) } ds \right\} \notag \\
						& = 1 -e^{-\lambda_1 t} + \sum_{k=2}^\infty \left(
						\prod_{j=1}^{k-1}
						\lambda_j \sum_{m=1}^k \frac{\left( 1-e^{-\lambda_m t} \right) }{
						\lambda_m \prod_{l=1,l\neq m}^k \left( \lambda_l -
						\lambda_m \right)} \right)
						\notag \\
						& = 1-e^{-\lambda_1 t} + \sum_{k=2}^\infty \sum_{m=1}^k
						\prod_{l=1,l\neq m}^k
						\frac{\lambda_l}{\lambda_l-\lambda_m}
						\left(1-e^{-\lambda_m t}\right) \notag \\
						& = 1 - e^{-\lambda_1 t} + \sum_{k=2}^\infty \left\{ 1-\sum_{m=1}^k
						\prod_{l=1,l\neq m}^k \frac{\lambda_l}{\lambda_l-\lambda_m}
						e^{-\lambda_m t} \right\} \notag \\
						& = \sum_{k=1}^\infty \left\{ 1-\sum_{m=1}^k
						\prod_{l=1,l\neq m}^k \frac{\lambda_l}{\lambda_l-\lambda_m}
						e^{-\lambda_m t} \right\}, \notag
						\notag
					\end{align}
					so that formula \eqref{valatt2} emerges.
					In the second-to-last step of \eqref{levis}, we applied formula (3.12) of
					\citet{sakhno} and, in the last step, we considered that,
					for $k=1$, the set of numbers
					$\{1\leq l \leq 1$, $l\neq m=1\}$, is empty and the
					\begin{equation}
						\prod_{l=1,l\neq m}^k \frac{\lambda_l}{\lambda_l -\lambda_m},
					\end{equation}
					is taken equal to 1 by convention.
				\end{proof}
			\end{thm}

			\begin{rem}
				As a check we can extract, from \eqref{valatt2}, the mean value in the linear case
				$\lambda_m = m \cdot \lambda$.
				Since
				\begin{align}
					& \sum_{m=1}^k \prod_{l=1, l\neq m}^k \frac{\lambda l}{\lambda l -
					\lambda m}
					e^{-\lambda m t}
					= \sum_{m=1}^k \frac{1 \dots (m-1)(m+1) \dots k}{(m-1)(m-2)
					\dots 1\cdot (-1)^{m-1}\cdot 1 \dots (k-m)} e^{-\lambda m t} \\
					& = - \sum_{m=1}^k \frac{k!}{m!(k-m)!} (-1)^m e^{-\lambda m t}
					= - \sum_{m=1}^k \binom{k}{m} (-1)^m e^{-\lambda m t}
					= - \left[ (1-e^{-\lambda t})^k -1 \right], \notag
				\end{align}
				we have that
				\begin{equation}
					1 - \sum_{m=1}^k \prod_{l=1,l\neq m}^k
					\frac{\lambda l}{\lambda l - \lambda m}
					e^{-\lambda m t} = (1-e^{-\lambda t})^k.
				\end{equation}
				From this we readily have that
				\begin{align}
					\mathbb{E} N(t) & = 1+\sum_{k=1}^\infty (1-e^{-\lambda t})^k
					= 1 + \frac{1}{1-(1-e^{-\lambda t})} -1 = e^{\lambda t}.
				\end{align}
			\end{rem}

			The aim of this section is to compose the process $\mathpzc{N}(t)$
			with the first-passage
			time $T_t = \inf (s \colon B(s)=t)$, where $B$ is a Brownian motion independent of
			$\mathpzc{N}(t)$.
			\begin{rem}
				\label{first}
				The probability density of $T_t= \inf \left\{ s\colon B(s)=t \right\}$,
				$t>0$, where $B(t)$ is a standard Brownian motion, namely
				\begin{equation}
					\text{Pr} \left\{ T_t \in ds \right\}/ds = q(t,s)
					=t\frac{e^{-\frac{t^2}{2s}}}{\sqrt{2\pi s^3}},
				\end{equation}
				is the solution to the Cauchy problem
				\begin{equation}
					\begin{cases}
						\frac{\partial^2}{\partial t^2} q(t,s) =
						2 \frac{\partial}{\partial s} q(t,s), \qquad t>0, \: s>0,\\
						q(0,s) = \delta(s),
					\end{cases}
				\end{equation}
				as a simple check shows.
			\end{rem}

		\subsection{Pure birth process stopped at $T_t$}

			\begin{thm}
				\label{theocomp}
				Let $\mathpzc{N} (t)$, $t>0$ be a classical nonlinear pure birth process and
				let $q(t,s)$, $s>0$, $t>0$, the law of $T_t$.
				The process $\hat{\mathpzc{N}}(t)=
				\mathpzc{N}(T_t)$, $t>0$, has the following distribution
				\begin{equation}
					\label{nlinear-comp}
					\hat{\mathpzc{p}}_k (t) =
					\begin{cases}
						\prod_{j=1}^{k-1} \lambda_j
						\sum_{m=1}^k \frac{ e^{-
						t \sqrt{2 \lambda_m}}}{\prod_{
						l=1,l \neq m}^k \left( \lambda_l -
						\lambda_m \right) }, & k > 1, \: t>0, \\
						e^{-t \sqrt{2 \lambda_1}}, & k=1, \: t>0,
					\end{cases}
				\end{equation}
				and mean value equal to
				\begin{equation}
					\label{meanT}
					\mathbb{E} \mathpzc{N}(T_t) = 1 + \sum_{k=1}^\infty \left(
					1 - \sum_{m=1}^k \prod_{l=1,l\neq m}^k
					\frac{\lambda_l}{\lambda_l - \lambda_m}
					e^{-t \sqrt{2\lambda_m}} \right).
				\end{equation}
				The distribution \eqref{nlinear-comp} is non-exploding under the condition
				that $\sum_k 1/\lambda_k =\infty$.

				\begin{proof}
					The state probabilities are derived by straight calculations and
					by resorting to the Laplace transform of $q(t,s)$ which reads
					\begin{equation}
						\int_0^\infty e^{- \gamma s} q(t,s) ds =
						\int_0^\infty e^{- \gamma s}
						\frac{t e^{- \frac{t^2}{2 s}}}{\sqrt{
						2 \pi s^3}} ds = e^{-t \sqrt{2 \gamma}}.
					\end{equation}
					We treat the case $k > 1$ as follows. The case $k=1$ is analogous.
					\begin{align}
						\label{b}
						\hat{\mathpzc{p}}_k (t) & = \text{Pr} \left\{
						\mathpzc{N} (T_t) = k \mid \mathpzc{N} (0)=1 \right\}
						= \int_0^\infty \mathpzc{p}_k (s) q(t,s) ds \\
						& = \int_0^\infty \prod_{j=1}^{k-1} \lambda_j
						\sum_{m=1}^k \frac{ e^{-
						\lambda_m s}}{\prod_{
						l=1,l \neq m}^k \left( \lambda_l -
						\lambda_m \right) } \frac{t e^{- \frac{t^2}{2 s}}}{\sqrt{
						2 \pi s^3}} ds
						= \prod_{j=1}^{k-1} \lambda_j
						\sum_{m=1}^k
						\frac{1}{\prod_{
						l=1,l \neq m}^k \left( \lambda_l -
						\lambda_m \right)} e^{- t \sqrt{2 \lambda_m}} . \nonumber
					\end{align}

					In view of Theorem \ref{ball}, we can evaluate the mean
					value
					\begin{align}
						\mathbb{E} \mathpzc{N}(T_t) & = \int_0^\infty \mathbb{E}
						\mathpzc{N}(s) \frac{te^{-\frac{t^2}{2s}}}{\sqrt{2\pi s^3}}
						ds
						= 1 + \sum_{k=1}^\infty \left(
						1 - \sum_{m=1}^k \prod_{l=1,l\neq m}^k
						\frac{\lambda_l}{\lambda_l - \lambda_m}
						e^{-t \sqrt{2\lambda_m}} \right).
					\end{align}
				\end{proof}
			\end{thm}

			In the linear case \eqref{meanT} can be written as
			\begin{equation}
				\mathbb{E} \mathpzc{N} (T_t) =
				\sum_{k=0}^\infty \sum_{m=0}^k (-1)^m e^{-t\sqrt{2\lambda m}}.
			\end{equation}
			On the other side, this sum diverges because
			\begin{equation}
				\mathbb{E} \mathpzc{N} (T_t) = \int_0^\infty e^{\lambda s}
				\text{Pr}\left\{ T_t \in ds \right\} = \infty.
			\end{equation}

			\begin{rem}
				Note that $\forall \, t$, $\hat{\mathpzc{p}}_k (t)$, $k \geq 1$ is a proper
				probability distribution because of the composition $\hat{\mathpzc{N}}(t)=
				\mathpzc{N}(T_t)$. The process can be appropriately interpreted
				by rewriting \eqref{b} as follows
				\begin{align}
					\hat{\mathpzc{p}}_k (t) & = \prod_{j=1}^{k-1} \lambda_j
					\sum_{m=1}^k
					\frac{1}{\prod_{
					l=1,l \neq m}^k \left( \lambda_l -
					\lambda_m \right) } \int_0^\infty e^{-\lambda_m s}
					\frac{t e^{- \frac{t^2}{2 s}}}{\sqrt{
					2 \pi s^3}} ds \\
					& = \int_0^\infty \prod_{j=1}^{k-1} \vartheta \lambda_j
					\sum_{m=1}^k
					\frac{1}{\prod_{
					l=1,l \neq m}^k \left( \vartheta \lambda_l -
					\vartheta \lambda_m \right) } e^{-\lambda_m \vartheta t^2}
					\frac{e^{- \frac{1}{2 \vartheta}}}{\sqrt{
					2 \pi \vartheta^3}} d \vartheta \nonumber .
				\end{align}
				The process $\hat{\mathpzc{N}}(t)$, $t>0$ can be viewed as a classical
				nonlinear pure birth process evaluated at time $t^2$ with random birth
				rates $\Theta \lambda_k$, $k \geq 1$, where $\Theta$ is an
				inverse Gaussian random variable
				with p.d.f.
				\begin{equation}
					f_\Theta (\vartheta) = \frac{e^{-\frac{1}{2 \vartheta}}}{\sqrt{
					2 \pi \vartheta^3}}, \qquad \vartheta \in \mathbb{R}^+ .
				\end{equation}
			\end{rem}

			The composition of $\mathpzc{N}(t)$, $t>0$, with $T_t$ leads to a second-order
			time derivative in the governing equations, as shown in the next theorem.

			\begin{thm}
				\label{solving}
				Let $\hat{\mathpzc{p}}_k (t)$, $t>0$, $k \geq 1$, be the distribution
				of the process $\hat{\mathpzc{N}} (t) = \mathpzc{N}(T_t)$,
				$t>0$, where $T_t$ is the first-passage time process of the
				standard Brownian motion, having transition density $q(t,s)$, $s>0$, $t>0$.
				The state probabilities
				$\hat{\mathpzc{p}}_k (t)$, $t>0$, $k \geq 1$, satisfy the following
				difference-differential equations
				\begin{equation}
					\label{eq-sec}
					\frac{d^2}{d t^2} \hat{\mathpzc{p}}_k (t)
					= 2 \left[ \lambda_k \hat{\mathpzc{p}}_k (t) - \lambda_{k-1}
					\hat{\mathpzc{p}}_{k-1} (t) \right], \qquad k \geq 1,
				\end{equation}
				where $\lambda_k$, $k \geq 1$ are the birth rates of the nonlinear
				classical birth process $\mathpzc{N} (t)$, $t>0$.
				\begin{proof}
					Since
					\begin{equation}
						\hat{\mathpzc{p}}_k (t) = \int_0^\infty
						\mathpzc{p}_k (s) q(t,s) ds,
					\end{equation}
					by taking the second-order derivative w.r.t.\ $t$,
					in view of Remark \ref{first}, we have that
					\begin{align}
						\label{microbiologicamente}
						\frac{d^2}{d t^2} \hat{\mathpzc{p}}_k (t)
						& = 2 \int_0^\infty \mathpzc{p}_k
						(s) \frac{\partial}{\partial s} q(t,s) ds
						= \left. 2 q(t,s) \mathpzc{p}_k (s)
						\right|_{s=0}^{s=\infty}
						- 2 \int_0^\infty \frac{d}{ds} \mathpzc{p}_k (s)
						q(t,s) ds \\
						& = -2 \int_0^\infty q(t,s) \left[ - \lambda_k
						\mathpzc{p}_k (s) + \lambda_{k-1} \mathpzc{p}_{k-1}
						(s) \right] ds
						= 2 \left[ \lambda_k \hat{\mathpzc{p}}_k (t)
						- \lambda_{k-1} \hat{\mathpzc{p}}_{k-1} (t)
						\right] . \nonumber
					 \end{align}
					 In \eqref{microbiologicamente}, we considered that
					 $\mathpzc{p}_k(0)=0$, for $k>1$.
				\end{proof}
			\end{thm}

			\begin{rem}
				In the linear case, some calculations suffice to show that
				\begin{equation}
					\hat{p}_k(t) = \text{Pr} \left\{ N(T_t) = k \right\}
					= \sum_{m=1}^k \binom{k-1}{m-1} (-1)^{m-1} e^{-t \sqrt{2\lambda m}},
					\qquad k\geq 1, \: t>0,
				\end{equation}
				and the state probabilities satisfy the equation
				\begin{equation}
					\frac{d^2}{dt^2}\hat{p}_k(t) = 2\lambda \hat{p}_k (t)
					-2\lambda (k-1) \hat{p}_{k-1} (t), \qquad k \geq 1.
				\end{equation}
			\end{rem}

			\subsubsection{Iterated compositions}

				\begin{thm}
					Let $\mathpzc{N} (t)$, $t>0$,
					be a classical nonlinear birth process.
					Let $T^1_t$, $T^2_t, \dots, T^n_t$, $n \in \mathbb{N}$,
					be first-passage times of $n$ independent
					standard Brownian motions.
					The process
					\begin{equation}
						\label{iterate}
						\tilde{\mathpzc{N}} (t) =
						\mathpzc{N} \biggl[ T^1_{T^2_{._{._{{T^n_t} }}}} \biggr],
						\qquad t>0,
					\end{equation}
					has the following distribution
					\begin{equation}
						\label{nlinear-mcomp}
						\tilde{\mathpzc{p}}_k (t) =
						\begin{cases}
							\prod_{j=1}^{k-1} \lambda_j
							\sum_{m=1}^k \frac{ e^{-t
							\lambda_m^{\frac{1}{
							2^n}} 2^{\left(1-\frac{1}{
							2^n} \right)}}}{\prod_{
							l=1,l \neq m}^k \left( \lambda_l -
							\lambda_m \right) }, & k > 1, \: t>0, \\
							e^{-t \lambda_1^{\frac{1}{
							2^n}} 2^{\left(1-\frac{1}{
							2^n} \right)}}, & k=1, \: t>0.
						\end{cases}
					\end{equation}
					\begin{proof}
						We start by proving the case $n=2$ since the
						case $n=1$ is already proved in Theorem \ref{theocomp}.
						We omit the details for the case $k=1$ and directly treat the
						case $k \geq 2$. We have that
						\begin{align}
							\label{compform}
							\int_0^\infty \hat{\mathpzc{p}}_k (s) q(t,s) ds
							& = \int_0^\infty
							\prod_{j=1}^{k-1} \lambda_j
							\sum_{m=1}^k \frac{ e^{-
							s \sqrt{2 \lambda_m}}}{\prod_{
							l=1,l \neq m}^k \left( \lambda_l -
							\lambda_m \right) } \frac{t e^{- \frac{t^2}{2 s}}}{\sqrt{
							2 \pi s^3}} ds \\
							& = \prod_{j=1}^{k-1} \lambda_j
							\sum_{m=1}^k
							\frac{1}{\prod_{
							l=1,l \neq m}^k \left( \lambda_l -
							\lambda_m \right) } \int_0^\infty e^{-s
							\sqrt{2 \lambda_m}}
							\frac{t e^{- \frac{t^2}{2 s}}}{\sqrt{
							2 \pi s^3}} ds \nonumber \\
							& = \prod_{j=1}^{k-1} \lambda_j
							\sum_{m=1}^k
							\frac{1}{\prod_{
							l=1,l \neq m}^k \left( \lambda_l -
							\lambda_m \right)}
							e^{- t \sqrt{2 \sqrt{2 \lambda_m}}} . \nonumber
						\end{align}
						It is now straightforward to generalise formula \eqref{compform}
						for $n$ compositions, as follows
						\begin{align}
							\label{beforeinfty}
							\tilde{\mathpzc{p}}_k (t) & =
							\prod_{j=1}^{k-1} \lambda_j
							\sum_{m=1}^k
							\frac{1}{\prod_{
							l=1,l \neq m}^k \left( \lambda_l -
							\lambda_m \right)}
							e^{- t \lambda_m^{\frac{1}{2^n}} 2^{\sum_{i=1}^n
							\frac{1}{2^i}}} \\
							& = \prod_{j=1}^{k-1} \lambda_j
							\sum_{m=1}^k
							\frac{1}{\prod_{
							l=1,l \neq m}^k \left( \lambda_l -
							\lambda_m \right)}
							e^{- t \lambda_m^{\frac{1}{2^n}} 2^{\left(
							1-\frac{1}{2^n} \right) }} . \nonumber
						\end{align}
					\end{proof}
				\end{thm}

				When $n \rightarrow \infty$, equation \eqref{nlinear-mcomp} becomes
				\begin{equation}
					\lim_{n \rightarrow \infty}\tilde{\mathpzc{p}}_k (t) =
					\begin{cases}
						e^{-2t}
						\prod_{j=1}^{k-1} \lambda_j
						\sum_{m=1}^k
						\frac{1}{\prod_{
						l=1,l \neq m}^k \left( \lambda_l -
						\lambda_m \right)} = 0, & k>1, \\
						e^{-2 t}, & k=1,
					\end{cases}
				\end{equation}
				because of formula (3.4), page 51 of \citet{chiang}.
				Therefore, the process \eqref{iterate} can either
				assume the state $k=1$ with
				probability $e^{-2t}$, or explode with probability $1-e^{-2t}$.

				\begin{thm}
					\label{solving-ite}
					Let $\tilde{\mathpzc{p}}_k (t)$, $t>0$, $k \geq 1$, be the distribution
					of the process
					\begin{equation}
						\label{iterate2}
						\tilde{\mathpzc{N}} (t) =
						\mathpzc{N} \biggl[ T^1_{T^2_{._{._{{T^n_t} }}}} \biggr],
						\qquad t>0.
					\end{equation}
					The state probabilities
					$\tilde{\mathpzc{p}}_k (t)$, $t>0$, $k \geq 1$, satisfy the following
					difference-differential equations
					\begin{equation}
						\label{olidata}
						\frac{d^{2^n}}{dt^{2^n}} \tilde{\mathpzc{p}}_k(t) =
						2^{2^n-1} \left\{ \lambda_k \tilde{\mathpzc{p}}_k(t) -
						\lambda_{k-1} \tilde{\mathpzc{p}}_{k-1}(t)\right\},
					\end{equation}
					where $\lambda_k$, $k \geq 1$, are the birth rates of the nonlinear
					classical birth process $\mathpzc{N} (t)$, $t>0$.
					\begin{proof}
						For $n=1$, equations \eqref{olidata} reduce to
						equations \eqref{eq-sec}. For $n=2$ we have that
						\begin{align}
							\frac{d^4}{dt^4} \tilde{\mathpzc{p}}_k(t)
							& = \int_0^\infty \int_0^\infty
							\hat{\mathpzc{p}}_k (w_1) q(w_2,w_1)
							\frac{\partial^4}{\partial t^4}
							q(t,w_2) dw_1 \, dw_2 \\
							& = 2^2 \int_0^\infty \int_0^\infty
							\hat{\mathpzc{p}}_k (w_1) \frac{\partial^2}{
							\partial w_2^2} q(w_2,w_1) q(t,w_2) dw_1 \, dw_2
							\notag \\
							& = 2^2 \int_0^\infty \int_0^\infty \hat{\mathpzc{p}}_k
							(w_1) \frac{\partial^2}{\partial w_2^2}
							q(w_2,w_1) q(t,w_2) dw_1 \, dw_2 \notag \\
							& = 2^3 \int_0^\infty \int_0^\infty
							\hat{\mathpzc{p}}_k (w_1) \frac{\partial}{\partial w_1}
							q(w_2,w_1) q(t,w_2) dw_1 \, dw_2 \notag \\
							& = -2^3 \int_0^\infty \int_0^\infty
							\frac{d}{dw_1} \hat{\mathpzc{p}}_k (w_1)
							q(w_2,w_1) q(t,w_2) dw_1 \, dw_2 \notag \\
							& = 2^3 \left\{ \lambda_k \tilde{\mathpzc{p}}_k(t)
							- \lambda_{k-1} \tilde{\mathpzc{p}}_{k-1} (t)\right\}.
							\notag
						\end{align}
						The above reasoning can be generalised, thus arriving
						at equation \eqref{olidata}.
					\end{proof}
				\end{thm}

		\subsection{Other compositions}

			In this part we present the distributions of the classical
			nonlinear birth process $\mathpzc{N} (t)$,
			$t>0$, stopped at various random time processes, namely
			the sojourn time $\Gamma_t$ of a standard Brownian motion,
			the sojourn time $\mathfrak{G}_t$ of a standard Brownian bridge
			and the stable subordinator $\mathpzc{S}^\alpha(t)$ of order $\alpha \in (0,1]$.

			We start first by considering the nonlinear birth process at time
			\begin{align}			
				\Gamma_t = \int_0^t I_{[0,\infty)} (B(s)) ds = \text{meas} \left\{
				s<t \colon B(s) > 0 \right\}.
			\end{align}
			The process $\mathpzc{N} (\Gamma_t)$, is a slowed
			down birth process. In the next theorem we provide its distribution.

			\begin{thm}
				We have that
				\begin{equation}
					\label{soj}
					\text{Pr} \left\{ \mathpzc{N}(\Gamma_t)=k \right\} =
					\begin{cases}
						\prod_{j=1}^{k-1} \lambda_j
						\sum_{m=1}^k \frac{ e^{-
						\frac{t}{2}\lambda_m} I_0\left( \frac{t}{2} \lambda_m
						\right) }{\prod_{
						l=1,l \neq m}^k \left( \lambda_l -
						\lambda_m \right) }, & k > 1, \: t>0, \\
						e^{-\frac{t}{2}\lambda_1} I_0\left( \frac{t}{2}\lambda_1
						\right) , & k=1, \: t>0,
					\end{cases}
				\end{equation}
				where
				\begin{equation}
					I_0(z) = \sum_{k=0}^\infty \left( \frac{z}{2} \right)^{2k}
					\frac{1}{(k!)^2},
				\end{equation}
				is the zero-order Bessel function with imaginary argument.
				\begin{proof}
					The derivation of \eqref{soj} is based on the evaluation
					of the following integral:
					\begin{equation}
						\int_0^t e^{-s\lambda_m} \frac{ds}{\pi \sqrt{s
						\left(t-s\right)}} = e^{-\frac{t}{2}\lambda_m}
						I_0 \left( \frac{t}{2}\lambda_m \right).
					\end{equation}
				\end{proof}
			\end{thm}

			\begin{rem}
				In view of the integral representation of the Bessel function
				\begin{equation}
					I_0(z) = \frac{1}{2\pi} \int_0^{2\pi} e^{z\cos \vartheta}
					d \vartheta,
				\end{equation}
				we can give the following alternative, interesting representation of
				\eqref{soj}.
				\begin{equation}
					\text{Pr} \left\{ \mathpzc{N}(\Gamma_t) = k \right\}
					= \frac{1}{2\pi} \int_0^{2\pi} \text{Pr} \left\{
					\mathpzc{N} \left( t\sin^2\frac{\vartheta}{2} \right) = k \right\}
					d \vartheta.
				\end{equation}
				In other words,
				\begin{equation}
					\mathpzc{N}(\Gamma_t) \overset{\text{i.d.}}{=}
					\mathpzc{N} \left( t\sin^2\frac{\Theta}{2} \right),
				\end{equation}
				where $\Theta$ is a random variable uniform in $[0,2\pi]$.
			\end{rem}

			\begin{thm}
				For the nonlinear birth process stopped at
				\begin{equation}
					\label{f3}
					\mathfrak{G}_t = \int_0^t I_{[0,\infty)} (\bar{B}(s)) ds,
				\end{equation}
				$\bar{B}(s)$, $s>0$, being a Brownian bridge, we have that
				\begin{equation}
					\label{brid}
					\text{Pr} \left\{ \mathpzc{N} (\mathfrak{G}_t) = k \right\} =
					\begin{cases}
						\frac{1}{\lambda_k t} \left\{ 1- \sum_{m=1}^k
						\prod_{l=1,l\neq m}^k \left(
						\frac{\lambda_l}{\lambda_l-\lambda_m}
						\right)
						e^{-\lambda_m t} \right\}, & k>1, \\
						\frac{1-e^{-\lambda_1 t}}{\lambda_1 t}, & k=1.
					\end{cases}
				\end{equation}
				\begin{proof}
					The calculation
					\begin{equation}
						\text{Pr} \left\{ \mathpzc{N} (\mathfrak{G}_t) = k \right\}
						= \prod_{j=1}^{k-1} \lambda_j \sum_{m=1}^k
						\frac{1}{\prod_{l=1,l\neq m}^k \left( \lambda_l
						-\lambda_m\right)} \int_0^t e^{-\lambda_m s} \text{Pr}
						\left\{ \mathfrak{G}_t \in ds \right\},
					\end{equation}
					is sufficient to arrive at result \eqref{brid}, once the
					well-known fact that \eqref{f3} is uniformly distributed in
					$[0,t]$ is considered.
				\end{proof}
			\end{thm}

			\begin{rem}
				For the linear birth process, the distribution \eqref{brid}
				takes a very simple form as the calculations below show.
				Since for $\lambda_k = \lambda \cdot k$, $\lambda>0$, we have that
				\begin{equation}
					\prod_{l=1,l\neq m}^k \frac{\lambda_l}{\lambda_l -\lambda_m} =
					\binom{k}{m} (-1)^{m-1},
				\end{equation}
				we can write that
				\begin{align}
					\label{bridlin}
					\text{Pr} \left\{ N ( \mathfrak{G}_t ) =k \right\}
					& = \frac{1}{\lambda k t} \left( 1- \sum_{m=1}^k \binom{k}{m}
					(-1)^{m-1} e^{-\lambda m t} \right) \\
					& = \frac{1}{\lambda k t} \sum_{m=0}^k \binom{k}{m} (-1)^m
					e^{-\lambda m t}
					= \frac{\left( 1-e^{-\lambda t} \right)^k}{\lambda k t},
					\qquad k \geq 1.
					\notag
				\end{align}
				The distribution \eqref{bridlin} is logarithmic with
				parameter $1-e^{-\lambda t}$.
				In the logarithmic distribution with parameter $0<q<1$, we have that
				\begin{equation}
					\mathbb{E} L = - \frac{q}{(1-q) \log (1-q)},
				\end{equation}
				\begin{equation}
					\mathbb{V}\text{ar} L = - \frac{q}{(1-q)^2\log(1-q)} \left[
					1+ \frac{q}{\log{(1-q)}}\right].
				\end{equation}
				In our case $q=1-e^{-\lambda t}$ so that
				\begin{equation}
					\mathbb{E} N(\mathfrak{G}_t) = \frac{e^{\lambda t}-1}{\lambda t},
				\end{equation}
				\begin{equation}
					\mathbb{V}\text{ar} N(\mathfrak{G}_t) = \frac{e^{\lambda t}
					(e^{\lambda t }-1)}{
					\lambda t} \left[ 1-\frac{1-e^{-\lambda t}}{\lambda t} \right] .
				\end{equation}
				For large values of $t$ we have that
				\begin{equation}
					\mathbb{E} N(\mathfrak{G}_t ) \sim \frac{e^{\lambda t}}{\lambda t} =
					\frac{\mathbb{E}N(t)}{\lambda t},
				\end{equation}
				\begin{equation}
					\mathbb{V}\text{ar} N(\mathfrak{G}_t) \sim \frac{e^{\lambda t}
					(e^{\lambda t}-1)}{\lambda t} =
					\frac{\mathbb{V}\text{ar}N(t)}{\lambda t}.
				\end{equation}
			\end{rem}

			\begin{thm}
				For the nonlinear birth process stopped at an $\alpha$-stable time
				$\mathpzc{S}^\alpha (t)$
				with distribution $q_\alpha (t,s)$ and Laplace transform $\int_0^\infty
				e^{-\mu s} q_\alpha (t,s) ds = e^{-t\mu^\alpha}$, we have that
				\begin{equation}
					\label{f4}
					\text{Pr} \left\{ \mathpzc{N} (\mathpzc{S}^\alpha (t)) = k \right\} =
					\begin{cases}
						\prod_{j=1}^{k-1} \lambda_j
						\sum_{m=1}^k \frac{
						e^{-t \lambda_m^\alpha}}{\prod_{
						l=1,l \neq m}^k \left( \lambda_l -
						\lambda_m \right) }, & k>1,\\
						e^{-t\lambda_1^\alpha}, & k=1.
					\end{cases}
				\end{equation}
				\begin{proof}
					The following calculation is sufficient to prove result \eqref{f4}:
					\begin{equation}
						\text{Pr} \left\{ \mathpzc{N} (\mathpzc{S}^\alpha (t))
						= k \right\} =
						\prod_{j=1}^{k-1} \lambda_j
						\sum_{m=1}^k \frac{1}{\prod_{
						l=1,l \neq m}^k \left( \lambda_l -
						\lambda_m \right) } \int_0^\infty e^{-\lambda_m s} q_\alpha(t,s)
						ds.
					\end{equation}
				\end{proof}
			\end{thm}

			\begin{rem}
				Formula \eqref{f4} can be further worked out as follows.
				\begin{align}
					\label{nlinear-sub}
					\text{Pr} \left\{ \mathpzc{N} (\mathpzc{S}^\alpha(t)) = k \right\} & =
					\prod_{j=1}^{k-1} \lambda_j
					\sum_{m=1}^k \frac{
					e^{-t \lambda_m^\alpha}}{\prod_{
					l=1,l \neq m}^k \left( \lambda_l -
					\lambda_m \right) } \\
					& \text{by exploiting the self-similarity of
					$\mathpzc{S}^\alpha(t)$} \notag \\
					& = \prod_{j=1}^{k-1} \lambda_j
					\sum_{m=1}^k \frac{1}{\prod_{
					l=1,l \neq m}^k \left( \lambda_l -
					\lambda_m \right) } \int_0^\infty e^{-\lambda_m s} t^{-\frac{1}{\alpha}}
					q_\alpha(1,t^{-\frac{1}{\alpha}} s) ds \notag \\
					& = \prod_{j=1}^{k-1} \lambda_j
					\sum_{m=1}^k \frac{1}{\prod_{
					l=1,l \neq m}^k \left( \lambda_l -
					\lambda_m \right) } \int_0^\infty e^{-\lambda_m t^{\frac{1}{\alpha}
					}\zeta} q_\alpha(1,\zeta) d \zeta. \notag
				\end{align}
				The last result implies the following representation:
				\begin{equation}
					\mathpzc{N}(\mathpzc{S}^\alpha(t)) \overset{\text{i.d.}}{=} \mathpzc{N}
					\left( t^{\frac{1}{\alpha}} Z \right),
				\end{equation}
				where $Z$ has distribution $q_\alpha (1,\zeta)$, $\zeta>0$.
			\end{rem}

			\begin{rem}
				If we assume $\alpha = 1/2^n$ in the first line of \eqref{nlinear-sub},
				and $s=t 2^{1-\frac{1}{2^n}}$ in \eqref{nlinear-mcomp}, the distribution
				\eqref{nlinear-sub} suggests the following unexpected relation:
				\begin{equation}
					\text{Pr} \left\{
					\mathpzc{N}\left(\mathpzc{S}^{\frac{1}{2^n}}(t) \right) = k
					\right\} = \text{Pr} \left\{
					\mathpzc{N} \biggl[ T^1_{T^2_{._{._{{T^n_s} }}}} \biggr]
					= k \right\}, \qquad k \geq 1.
				\end{equation}
			\end{rem}

			\begin{rem}
				Many other compositions can be envisaged and in some cases they
				provide curious results. For example, we consider the standard
				Cauchy process $C(t)$, with law $h(t,s)$, $t>0$, $s \in \mathbb{R}$,
				satisfying the Laplace equation
				\begin{equation}
					\frac{\partial^2 h}{\partial t^2} + \frac{\partial^2 h}{\partial s^2}
					= 0.
				\end{equation}
				We can show that $\mathpzc{N}(|C(t)|)$, $t>0$, is a birth process whose state
				probabilities $\mathpzc{p}_k^*(t)$, $t>0$, satisfy the difference-differential
				equations
				\begin{equation}
					\frac{d^2}{dt^2} \mathpzc{p}_k^* (t) =
					- \lambda_k^2 \mathpzc{p}_k^* (t) + \lambda_{k-1} \left(
					\lambda_k + \lambda_{k-1} \right) \mathpzc{p}_{k-1}^* (t) - \lambda_{k-1}
					\lambda_{k-2} \mathpzc{p}_{k-2}^* (t) .
				\end{equation}
			\end{rem}

	\section{Subordinated fractional birth processes}

		In a previous work of us (see \citet{pol}) we constructed and analysed a fractional
		(possibly nonlinear) pure birth process $\mathpzc{N}^\nu (t)$, $t>0$, $\nu \in (0,1]$
		by exchanging the integer-order time derivative with the Dzhrbashyan--Caputo
		fractional derivative in the difference-differential equation \eqref{fra}
		governing the state probabilities.
		We recall that the Dzhrbashyan--Caputo derivative has the form, for $0<\nu\leq 1$
		\begin{equation}
			\label{caputo}
			\frac{d^\nu}{dt^\nu} f(t) =
			\begin{cases}
				\frac{1}{\Gamma(1-\nu)} \int_0^t \frac{f'(s)}{(t-s)^\nu} ds, & 0<\nu<1, \\
				f'(t), & \nu=1.
			\end{cases}
		\end{equation}
		In this section we examine properties
		of the subordinated processes $\mathpzc{N}^\nu (T_t)$, $t>0$, $\mathpzc{N}^\nu
		(T_{2\beta}(t))$, and
		$\mathpzc{N}^\nu (\mathpzc{S}^\alpha(t))$, $t>0$, $\nu,\alpha,\beta \in (0,1]$,
		bringing to the fore some interesting relations and discussing the interpretation for
		the results obtained.

		\subsection{Preliminaries}

			The state probabilities $\mathpzc{p}_k^\nu (t)
			= \text{Pr} \left\{ \mathpzc{N}^\nu (t) = k \right\}$, $k \geq 1$ of the fractional
			pure birth process have the following form
			\begin{equation}
				\label{nlinearnu}
				\mathpzc{p}_k^\nu (t) =
				\begin{cases}
					\prod_{j=1}^{k-1} \lambda_j
					\sum_{m=1}^k \frac{ E_{\nu,1} (-
					\lambda_m t^\nu)}{\prod_{
					l=1,l \neq m}^k \left( \lambda_l -
					\lambda_m \right) }, & k > 1, \: t>0, \\
					E_{\nu,1} (-\lambda_1 t^\nu), & k=1, \: t>0,
				\end{cases}
			\end{equation}
			where $E_{\nu,1} (- \zeta t^\nu)$ is the Mittag--Leffler function defined as
			\begin{equation}
				\label{mittag}
				E_{\nu,1} \left( - \zeta t^\nu \right) = \sum_{h=0}^\infty
				\frac{(- \zeta t^\nu)^h}{\Gamma \left( \nu h+1 \right)}, \qquad
				\zeta \in \mathbb{R}, \: \nu > 0,
			\end{equation}
			and with Laplace transform
			\begin{equation}
				\label{mitlap}
				\int_0^\infty e^{- z t} E_{\nu, 1} \left( - \zeta t^\nu \right) dt =
				\frac{z^{\nu -1}}{z^\nu + \zeta}, \qquad \nu >0 .
			\end{equation}
			A useful integral representation for $E_{\nu,1} (- \zeta t^\nu)$ reads
			\begin{equation}
				\label{mitta-int}
				E_{\nu,1} (- \zeta t^\nu) = \frac{\sin \nu \pi}{\pi}
				\int_0^\infty \frac{r^{\nu-1} e^{-r \zeta^{\frac{1}{\nu}} t}}{
				r^{2 \nu} + 2 r^\nu \cos \nu \pi + 1} dr, \qquad \nu \in (0,1) .
			\end{equation}

			In a previous work (see \citet{pol}) we proved a useful
			subordination representation for the fractional pure birth process
			\eqref{nlinearnu}. This can be viewed as a classical birth process
			stopped at a random time $T_{2 \nu}(t)$ possessing density function
			coinciding with the folded solution to the fractional
			diffusion equation
			\begin{equation}
				\label{diffusion}
				\begin{cases}
					\frac{\partial^{2 \nu} g}{\partial t^{2 \nu}}
					= \frac{\partial^2 g}{\partial x^2},
					& 0 < \nu \leq 1, \\
					g \left( x, 0 \right) = \delta \left( x \right),
				\end{cases}
			\end{equation}
			with the additional condition $g_t \left(x,0 \right) = 0$ for $1/2 < \nu \leq 1$.
			In other words $\mathpzc{N}^\nu (t) = \mathpzc{N} (T_{2 \nu} (t) )$, $t>0$.
			It can be shown that
			$f_{T_{2\nu}}(s,t) = \text{Pr}\left\{ T_{2\nu}(t) \in ds \right\}$ is
			also a solution to
			\begin{equation}
				\frac{\partial^\nu f}{\partial t^\nu} = - \frac{\partial f}{\partial s}
			\end{equation}
			(see \citet{sakhno}).

			\begin{thm}
				The fractional nonlinear pure birth process is a renewal process with
				intermediate waiting times $T^\nu_k$ with law
				\begin{equation}
					\text{Pr} \left\{ T^\nu_k \in ds \right\} =
					\lambda_k s^{\nu-1} E_{\nu,\nu} (-\lambda_k s^\nu) ds,
					\qquad k\geq 1, \: s>0,
				\end{equation}
				where $T_k^\nu$ is the random time separating the $k$th and $(k+1)$th
				birth.

				\begin{proof}
					We prove this result by induction. Denoting $Z^\nu_k
					=T_1^\nu+\dots+T_k^\nu$, we can certainly write that
					\begin{align}
						\label{tex}
						\text{Pr}\left\{ T_1^\nu+\dots+T_k^\nu\in dt \right\}
						= \int_0^t \text{Pr} \left\{ T_k^\nu \in d(t-s) \right\}
						\text{Pr} \left\{ T_1^\nu+\dots+T^\nu_{k-1}\in ds \right\},
					\end{align}
					where $\text{Pr} \left\{ T_1^\nu+\dots+T_{k-1}^\nu\in ds \right\}
					/ ds = \frac{d}{ds} \text{Pr}\left\{ \mathpzc{N}^\nu(s)\geq k\right\}$.
					By resorting to Laplace transforms, from \eqref{tex}, we obtain
					that
					\begin{align}
						\label{capcap}
						\int_0^\infty e^{-\mu t} \text{Pr}
						\left\{ T_1^\nu+\dots+T_{k}^\nu\in dt \right\}
						& = \int_0^\infty e^{-\mu t} dt \int_0^t
						\text{Pr} \left\{ T_k^\nu\in d(t-s) \right\}
						\text{Pr}\left\{ T_1^\nu+\dots+T_{k-1}^\nu\in ds \right\}
						\\
						& = \int_0^\infty \text{Pr} \left\{ T_1^\nu +\dots+T_{k-1}^\nu
						\in ds \right\} \int_s^\infty e^{-\mu t}
						\text{Pr} \left\{ T_k^\nu \in d(t-s) \right\} \notag \\
						& = \int_0^\infty e^{-\mu s} \text{Pr} \left\{
						T_1^\nu+\dots+T_{k-1}^\nu\in ds \right\} \int_0^\infty
						e^{-\mu y} \text{Pr} \left\{ T_k^\nu \in dy \right\} \notag \\
						& = \prod_{j=1}^k \int_0^\infty e^{-\mu s}
						\text{Pr} \left\{ T_j^\nu\in ds \right\} =
						\prod_{j=1}^k \frac{\lambda_j}{\mu^\nu+\lambda_j}. \notag
					\end{align}
					We observe that
					\begin{align}
						\text{Pr}\left\{ T_1^\nu\in ds \right\}/ds & =
						\frac{d}{ds} \text{Pr}\left\{ \mathpzc{N}^\nu(s)\geq 2\right\}
						= -\frac{d}{ds} E_{\nu,1}(-\lambda_1s^\nu)
						= \lambda_1 s^{\nu-1} E_{\nu,\nu} (-\lambda_1 s^\nu),
					\end{align}
					and that
					\begin{align}
						\label{ebal}
						\text{Pr} \left\{ T_1^\nu +T_2^\nu \in ds \right\}
						& = \frac{d}{ds} \left[ 1-\text{Pr} \left\{
						\mathpzc{N}^\nu(s)=1 \right\} - \text{Pr} \left\{
						\mathpzc{N}^\nu(s)=2 \right\}\right] \\
						& = \lambda_1 s^{\nu-1} E_{\nu,\nu} (-\lambda_1 s^\nu)
						+ \lambda_1 \left[ \lambda_1 \frac{E_{\nu,\nu}(-\lambda_1
						s^\nu)}{\lambda_2-\lambda_1} + \lambda_2 \frac{
						E_{\nu,\nu}(-\lambda_2 s^\nu)}{\lambda_1-\lambda_2}\right]
						s^{\nu-1} \notag \\
						& = \frac{\lambda_1\lambda_2}{\lambda_2-\lambda_1} s^{\nu-1}
						\left[ E_{\nu,\nu}(-\lambda_1 s^\nu) - E_{\nu,\nu}
						(-\lambda_2 s^\nu)\right].\notag
					\end{align}
					For $k=2$, relation \eqref{capcap} simplifies to
					\begin{align}
						\int_0^\infty e^{-\mu t} \text{Pr} \left\{
						T_1^\nu+T_2^\nu\in dt \right\} & = \frac{\lambda_1}{\mu^\nu
						+\lambda_1} \frac{\lambda_2}{\mu^\nu+\lambda_2},
					\end{align}
					and this coincides with the Laplace transform of \eqref{ebal}.
				\end{proof}
			\end{thm}

			\begin{thm}
				The mean value $\mathbb{E} \mathpzc{N}^\nu (t)$, for the fractional
				nonlinear pure birth process has the form:
				\begin{equation}
					\label{meanmean}
					\mathbb{E} \mathpzc{N}^\nu (t) = 1 + \sum_{k=1}^\infty
					\left\{ 1- \sum_{m=1}^k \prod_{l=1,l\neq m}^k \frac{\lambda_l}{
					\lambda_l - \lambda_m} E_{\nu,1}(-\lambda_m t^\nu) \right\}.
				\end{equation}
				\begin{proof}
					In light of the subordination relation $\mathpzc{N}^\nu(t)
					= \mathpzc{N}(T_{2\nu}(t))$, and of result \eqref{valatt2},
					we can write that
					\begin{align}
						\mathbb{E} \mathpzc{N}^\nu (t) & =
						\int_0^\infty \mathbb{E} \mathpzc{N}(s)
						\text{Pr} \left\{ T_{2\nu} (t) \in ds \right\} \\
						& = 1 + \sum_{k=1}^\infty \left\{ 1- \sum_{m=1}^k
						\prod_{l=1,l\neq m}^k \frac{\lambda_l}{\lambda_l - \lambda_m}
						\int_0^\infty e^{-\lambda_m s} \text{Pr} \left\{
						T_{2\nu}(t) \in ds \right\} \right\} \notag \\
						& = 1 + \sum_{k=1}^\infty \left\{ 1 - \sum_{m=1}^k
						\prod_{l=1,l\neq m}^k \frac{\lambda_l}{\lambda_l-\lambda_m}
						E_{\nu,1}(-\lambda_m t^\nu) \right\}. \notag
					\end{align}
				\end{proof}
			\end{thm}
			In the previous steps we assumed that
			\begin{equation}
				\label{ackard}
				\int_0^\infty e^{-\lambda_m s} f_{T_{2\nu}} (t,s) ds
				= \int_0^\infty e^{-\lambda_m s} \text{Pr} \left\{
				T_{2\nu}(t) \in ds \right\}
				= E_{\nu,1} (-\lambda_m t^\nu).
			\end{equation}
			We give here some details of this result. The density $f_{T_{2\nu}}(
			z,s)$, $z>0$, $s>0$, is obtained by folding the solution
			of the fractional diffusion equation
			\begin{equation}
				\frac{\partial^{2\nu} u}{\partial z^{2\nu}} =
				\frac{\partial^2 u}{\partial s^2},
			\end{equation}
			which reads
			\begin{equation}
				u(z,s) = \frac{1}{z^\nu} W_{-\nu,1-\nu} \left(
				-\frac{s}{z^\nu} \right), \qquad s>0,\: z>0,
			\end{equation}
			where $W_{-\nu, 1-\nu}(-\xi)$ is a Wright function defined as
			\begin{equation}
				W_{-\nu, 1-\nu}(-\xi) = \sum_{r=0}^\infty \frac{(- \xi)^r}{r!
				\Gamma \left( 1- \nu (r+1) \right)} .
			\end{equation}

			Therefore
			\begin{align}
				\frac{1}{z^\nu} \int_0^\infty e^{-\lambda_m s}
				W_{-\nu,1-\nu} \left( -\frac{s}{z^\nu} \right)
				ds
				& = \frac{1}{z^\nu} \int_0^\infty e^{-\lambda_m s}
				\sum_{k=0}^\infty \left( -\frac{s}{z^\nu} \right)^k
				\frac{1}{k! \Gamma (-\nu k+1 -\nu )} ds \\
				& = \frac{1}{z^\nu} \sum_{k=0}^\infty
				\frac{(-1)^k}{\lambda_m^{k+1}\Gamma(-\nu k +1 -\nu)}
				\frac{1}{(z^\nu)^k}
				= \frac{1}{\lambda_m z^\nu} E_{-\nu,1-\nu}
				\left( -\frac{1}{\lambda_m z^\nu} \right) \notag \\
				& \text{(by formula (5.1) page 1825, \citet{orsbeg2})} \notag \\
				& = E_{\nu,1} (-\lambda_m z^\nu). \notag
			\end{align}

			\begin{rem}
				We can extract, from \eqref{meanmean}, the mean value of the fractional
				linear birth process obtained in \citet{pol}, formula (3.42), as follows.
				By considering that $\lambda_m = \lambda \cdot m$, formula \eqref{meanmean}
				becomes
				\begin{align}
					\label{becomes}
					\mathbb{E} N^\nu (t) & = 1 + \sum_{k=1}^\infty
					\left\{ 1- \sum_{m=1}^k (-1)^{m-1} \frac{k!}{m!(k-m)!}
					E_{\nu,1}(-\lambda m t^\nu) \right\} \\
					& = 1 + \sum_{k=1}^\infty \left\{ 1+ \sum_{m=1}^k
					(-1)^m \binom{k}{m} E_{\nu,1}(-\lambda m t^\nu) \right\} \notag \\
					& = 1 + \sum_{k=1}^\infty \sum_{m=0}^k \binom{k}{m} (-1)^m
					E_{\nu,1}(-\lambda m t^\nu)
					= \sum_{k=0}^\infty \sum_{m=0}^k \binom{k}{m} (-1)^m
					E_{\nu,1} (-\lambda m t^\nu). \notag
				\end{align}
				In order to obtain the desired result we pass to Laplace transforms
				and extract from \eqref{becomes} that
				\begin{align}
					\int_0^\infty e^{-\mu t} \mathbb{E} N^\nu(t) dt & =
					\sum_{k=0}^\infty \sum_{m=0}^k \binom{k}{m} (-1)^k
					\int_0^\infty e^{-\mu t} E_{\nu,1}(-\lambda m t^\nu) dt \\
					& = \sum_{k=0}^\infty \sum_{m=0}^k \binom{k}{m} (-1)^k
					\frac{\mu^{\nu-1}}{\mu^\nu + \lambda m}
					= \mu^{\nu-1} \int_0^\infty \sum_{k=0}^\infty \sum_{m=0}^k
					\binom{k}{m} (-1)^k e^{-w\mu^\nu} e^{-w\lambda m} dw \notag \\
					& = \mu^{\nu-1} \int_0^\infty \sum_{k=0}^\infty e^{-w\mu^\nu}
					\left( 1-e^{-w\lambda} \right)^k dw
					= \mu^{\nu-1} \int_0^\infty \frac{e^{-w\mu^\nu}}{1-\left(
					1- e^{-w\lambda} \right)} dw \notag \\
					&= \mu^{\nu-1} \int_0^\infty e^{-w \mu^\nu+w\lambda} dw
					= \frac{\mu^{\nu-1}}{\mu^\nu-\lambda}. \notag
				\end{align}
				By inverting the Laplace transform above, we can conclude that
				\begin{equation}
					\mathbb{E} N^\nu(t) = E_{\nu,1}(\lambda t^\nu),
				\end{equation}
				thus confirming our previous result.
			\end{rem}

			Here we remark that another interpretation in terms of random birth rates
			can be highlighted. If we write
			\begin{align}
				\label{a}
				\mathpzc{p}_k^\nu (t) & = \text{Pr} \left\{ \mathpzc{N}^\nu (t) =k
				\mid \mathpzc{N}^\nu (0) = 1 \right\} =
				\int_0^\infty \mathpzc{p}_k (s) \text{Pr} \left\{ T_{2 \nu}
				(t) \in ds \right\} \\
				& = \prod_{j=1}^{k-1} \lambda_j
				\sum_{m=1}^k \frac{1}{\prod_{
				l=1,l \neq m}^k \left( \lambda_l -
				\lambda_m \right) } \int_0^\infty e^{- \lambda_m s}
				\text{Pr} \left\{ T_{2 \nu} (t) \in ds \right\} \nonumber \\
				& = \prod_{j=1}^{k-1} \lambda_j
				\sum_{m=1}^k \frac{1}{\prod_{
				l=1,l \neq m}^k \left( \lambda_l -
				\lambda_m \right) } \int_0^\infty e^{- \lambda_m s} \frac{1}{t^\nu}
				W_{-\nu,1-\nu} \left( - \frac{s}{t^\nu} \right) ds \nonumber \\
				& = \prod_{j=1}^{k-1} \lambda_j
				\sum_{m=1}^k \frac{1}{\prod_{l=1,l \neq m}^k \left( \lambda_l -
				\lambda_m \right) } \int_0^\infty e^{- \lambda_m \xi t^\nu}
				W_{-\nu, 1-\nu} (- \xi) d \xi \nonumber \\
				& = \int_0^\infty W_{- \nu, 1-\nu} (- \xi) \text{Pr}
				\left\{ \mathpzc{N}_{\xi} (t^\nu) = k \mid \mathpzc{N}_{\xi} (0) = 1
				\right\} d \xi, \nonumber
			\end{align}
			we have that a fractional nonlinear pure birth process can be considered as a
			classical nonlinear pure birth process evaluated at a rescaled time $t^\nu$
			and with random rates $\lambda_k \Xi$, $k \geq 1$, where $\Xi$ is a
			random variable with density function
			\begin{equation}
				f_{\Xi} (\xi) = W_{-\nu,1-\nu} (-\xi), \qquad \xi \in \mathbb{R}^+.
			\end{equation}
			From \eqref{a}, the following interpretation also holds:
			\begin{equation}
				\mathpzc{N}^\nu(t) \overset{\text{i.d.}}{=} \mathpzc{N} (\Xi t^\nu).
			\end{equation}

			Note also that, from \eqref{mitta-int} and \eqref{a}, we have that
			\begin{equation}
				\label{mmm}
				E_{\nu,1} (-zt^\nu) = \int_0^\infty e^{-rzt^\nu} W_{-\nu,1-\nu}(-r) dr =
				\frac{\sin \nu \pi}{\pi} \int_0^\infty e^{-rz^{\frac{1}{\nu}}t}
				\frac{r^{\nu-1}}{r^{2\nu}+2r^\nu\cos \nu \pi+1} dr,
			\end{equation}
			which illustrates an interesting relation between the Wright function
			and the law of $\mathcal{W}_1$ (see \eqref{wallet}).
			Equation \eqref{mmm} can be derived directly as follows.
			\begin{align}
				& \int_0^\infty e^{-\gamma x} \frac{1}{\lambda t^\nu}
				W_{-\nu,1-\nu} \left( -\frac{x}{\lambda t^\nu} \right) dx
				= \frac{1}{\lambda t^\nu} \sum_{m=0}^\infty
				\frac{(-1)^m}{m! \Gamma\left( -\frac{\nu}{m} +1 -\nu \right)}
				\frac{1}{\left( \lambda t^\nu
				\right)^m} \int_0^\infty e^{-\gamma x} x^m dx \\
				& = \frac{1}{\lambda t^\nu} \sum_{m=0}^\infty
				\frac{(-1)^m}{\Gamma\left( -\frac{\nu}{m} +1 -\nu \right)}
				\frac{1}{\gamma \left( \gamma \lambda t^\nu \right)^m }
				= \frac{1}{\gamma \lambda t^\nu} E_{-\nu,1-\nu}
				\left( -\frac{1}{\gamma \lambda t^\nu} \right) \notag \\
				& (\text{by formula (5.1), page 1825, \citet{orsbeg2}}) \notag \\
				& = E_{\nu,1} (-\gamma \lambda t^\nu)
				= \frac{\sin \nu \pi}{\pi} \int_0^\infty
				\frac{e^{-\left( \gamma \lambda \right)^{\frac{1}{\nu} } tr }r^{\nu-1} }{
				\left( r^\nu + \cos \nu\pi \right)^2 + \sin^2 \nu\pi } dr. \notag
			\end{align}
			This yields \eqref{mmm} for $z=\gamma \lambda$.

			For further details on fractional pure birth process the reader can refer to \citet{pol}
			while Mittag--Leffler functions are extensively analysed in \citet{kilbas}.

		\subsection{Fractional pure birth process stopped at $T_t$}

			We consider here the composition of a fractional nonlinear pure birth process,
			denoted as $\mathpzc{N}^\nu (t)$, $t>0$, $\nu \in (0,1]$ with the
			first-passage time $T_t$ of a standard Brownian motion.
			In the following theorem we derive an interesting integral representation
			for the state probabilities $\hat{\mathpzc{p}}_k^\nu (t)$, $t>0$, $k \geq 1$, of
			$\hat{\mathpzc{N}}^\nu (t) = \mathpzc{N}^\nu (T_t)$, $t>0$, $\nu \in (0,1)$.
			\begin{thm}
				\label{teor}
				Let $\mathpzc{N}^\nu (t)$, $t>0$, $\nu \in (0,1)$, be a fractional
				nonlinear pure birth process and $T_t$ be the first-passage time
				process of the standard Brownian motion with distribution $q(t,s)$.
				The state probabilities $\hat{\mathpzc{p}}_k^\nu (t) =
				\text{Pr} \left\{ \mathpzc{N}^\nu (T_t) = k \mid \mathpzc{N}^\nu
				(0) = 1 \right\}$ possess the following integral form
				\begin{equation}
					\label{frco}
					\hat{\mathpzc{p}}_k^\nu (t) =
					\begin{cases}
						\prod_{j=1}^{k-1} \lambda_j
						\sum_{m=1}^k
						\frac{1}{\prod_{
						l=1,l \neq m}^k \left( \lambda_l -
						\lambda_m \right) }
						\frac{1}{i \pi} \int_0^\infty
						\frac{E_{2 \nu,1} \left( -x^{2 \nu} e^{i \pi \nu} \right)
						- E_{2 \nu,1} \left( -x^{2 \nu} e^{-i \pi \nu}
						\right)}{x+ t \sqrt{2 \lambda_m^{\frac{1}{\nu}}}} dx,
						& k > 1, \\
						\frac{1}{i \pi} \int_0^\infty
						\frac{E_{2 \nu,1} \left( -x^{2 \nu} e^{i \pi \nu} \right)
						- E_{2 \nu,1} \left( -x^{2 \nu} e^{-i \pi \nu}
						\right)}{x+ t \sqrt{2 \lambda_1^{\frac{1}{\nu}}}} dx, & k=1.
					\end{cases}
				\end{equation}
				\begin{proof}
					It is sufficient to prove \eqref{frco} in the case $k>1$, since the
					case $k=1$ is analogous. We have
					\begin{align}
						\label{inte1}
						\hat{\mathpzc{p}}_k^\nu (t) & = \int_0^\infty
						\mathpzc{p}_k^\nu (s)
						q(t,s) ds
						= \prod_{j=1}^{k-1} \lambda_j
						\sum_{m=1}^k
						\frac{1}{\prod_{
						l=1,l \neq m}^k \left( \lambda_l -
						\lambda_m \right) } \int_0^\infty E_{\nu,1} ( - \lambda_m s^\nu )
						\frac{t e^{- \frac{t^2}{2s}}}{\sqrt{2 \pi s^3}} ds.
					\end{align}
					In order to prove \eqref{frco}, by taking into
					consideration formula \eqref{mitta-int}, we do the following calculations
					\begin{align}
						& \int_0^\infty E_{\nu,1} (- \lambda_m s^\nu)
						\frac{t e^{- \frac{t^2}{2s}}}{\sqrt{2 \pi s^3}} ds
						= \int_0^\infty \frac{ t e^{-\frac{t^2}{2s}}}{\sqrt{2 \pi s^3}}
						\int_0^\infty \frac{\sin \nu \pi}{\pi} \frac{r^{\nu-1}
						e^{-r \lambda_m^{\frac{1}{\nu}} s}}{
						r^{2 \nu} + 2 r^\nu \cos \nu \pi + 1} dr \, ds \\
						& = \frac{\sin \nu \pi}{\pi} \int_0^\infty
						\frac{r^{\nu-1}}{r^{2 \nu} + 2r^\nu \cos \nu \pi +1}
						\int_0^\infty e^{-r \lambda_m^{\frac{1}{\nu}} s }
						\frac{t e^-\frac{t^2}{2s}}{\sqrt{2 \pi s^3}} ds \, dr
						\nonumber \\
						& = \frac{\sin \nu \pi}{\pi} \int_0^\infty
						\frac{r^{\nu -1}}{r^{2 \nu} + 2 r^\nu \cos \nu \pi +1}
						e^{- t \sqrt{2 r \lambda_m^{\frac{1}{\nu}}}} dr
						= \frac{\sin \nu \pi}{\pi} \int_0^\infty
						\frac{r^{\nu -1} e^{-t \sqrt{2 r \lambda_m^{\frac{1}{\nu}}}}}{
						\left( r^\nu + e^{i \pi \nu} \right) \left( r^\nu
						+ e^{- i \pi \nu} \right)} dr \nonumber \\
						& = \frac{1}{2 i \pi} \int_0^\infty \left( \frac{1}{r^\nu +
						e^{i \pi \nu}} - \frac{1}{r^\nu + e^{-i \pi \nu}} \right)
						r^{\nu -1} e^{-t \sqrt{2 r \lambda_m^{\frac{1}{\nu}}}} dr
						\nonumber \\
						& = \frac{1}{i \pi} \int_0^\infty \left( \frac{1}{
						w^{2 \nu} + e^{i \pi \nu}} - \frac{1}{w^{2 \nu} +
						e^{- i \pi \nu}} \right) w^{2 \nu -1}
						e^{-t w \sqrt{2 \lambda_m^{\frac{1}{\nu}}}} dw \nonumber \\
						& = \frac{1}{i \pi} \int_0^\infty \left[ \frac{w^{2 \nu-1}}{
						w^{2 \nu} + e^{i \pi \nu}} - \frac{w^{2 \nu -1}}{
						w^{2 \nu} + e^{- i \pi \nu}} \right] e^{-tw \sqrt{2
						\lambda_m^{\frac{1}{\nu}}}} dw \nonumber .
					\end{align}
					By using the Laplace transform \eqref{mitlap} we obtain
					\begin{align}
						\label{inte}
						& \int_0^\infty E_{\nu,1} (- \lambda_m s^\nu)
						\frac{t e^{- \frac{t^2}{2s}}}{\sqrt{2 \pi s^3}} ds \\
						& = \frac{1}{i \pi} \int_0^\infty \left[
						\int_0^\infty e^{-wx} E_{2 \nu,1} (-x^{2 \nu}e^{i \pi \nu}) dx
						- \int_0^\infty e^{-wx} E_{2 \nu,1} ( -x^{2 \nu} e^{-i \pi \nu})
						dx \right] e^{-t w \sqrt{2 \lambda_m^{\frac{1}{\nu}}}} dw
						\nonumber \\
						& = \frac{1}{i \pi} \int_0^\infty \frac{ E_{2 \nu,1} \left(
						- x^{2 \nu} e^{i \pi \nu} \right) - E_{2 \nu,1} \left(
						- x^{2 \nu} e^{-i \pi \nu} \right)}{x+t \sqrt{2
						\lambda_m^{\frac{1}{
						\nu}}}} dx . \nonumber
					\end{align}
					Formula \eqref{frco} is then proved by combining
					\eqref{inte1} and \eqref{inte}.
				\end{proof}
			\end{thm}
			\begin{rem}
				If $\lambda_k=k \lambda$, $k \geq 1$,
				the state probabilities $\hat{p}_k^\nu (t) = \text{Pr}
				\left\{ N^\nu (T_t) = k \mid N^\nu (0) = 1 \right\}$
				of a fractional linear pure birth process
				stopped at $T_t$ read
				\begin{equation}
					\label{frcolin}
					\hat{p}_k^\nu (t) =
					\begin{cases}
						\sum_{m=1}^k
						\binom{k-1}{m-1} (-1)^{m-1} \frac{1}{i \pi} \int_0^\infty
						\frac{E_{2 \nu,1} \left( -x^{2 \nu} e^{i \pi \nu} \right)
						- E_{2 \nu,1} \left( -x^{2 \nu} e^{-i \pi \nu}
						\right)}{x+ t \sqrt{2 \lambda^{\frac{1}{\nu}} m }} dx,
						& k > 1, \\
						\frac{1}{i \pi} \int_0^\infty
						\frac{E_{2 \nu,1} \left( -x^{2 \nu} e^{i \pi \nu} \right)
						- E_{2 \nu,1} \left( -x^{2 \nu} e^{-i \pi \nu}
						\right)}{x+ t \sqrt{2 \lambda^{\frac{1}{\nu}}}} dx, & k=1.
					\end{cases}
				\end{equation}
				This result can be obtained by means of methods similar to
				those of Theorem \ref{teor}.
			\end{rem}
			\begin{rem}
				\label{interpretation}
				By considering formula \eqref{inte1} and the representation
				\eqref{mitta-int} we can give an interesting interpretation of the process
				$\hat{\mathpzc{N}}^{1/2} (t) = \mathpzc{N}^{1/2} (T_t)$,
				$t>0$, as follows (again, we treat the case $k \geq 1$ since
				the case $k=1$ is analogous)
				\begin{align}
					\label{inte2}
					\hat{\mathpzc{p}}_k^\nu (t) & = \int_0^\infty \mathpzc{p}_k^\nu (s)
					q(t,s) ds \\
					& = \prod_{j=1}^{k-1} \lambda_j
					\sum_{m=1}^k
					\frac{1}{\prod_{
					l=1,l \neq m}^k \left( \lambda_l -
					\lambda_m \right) } \int_0^\infty E_{\nu,1} ( - \lambda_m s^\nu )
					\frac{t e^{- \frac{t^2}{2s}}}{\sqrt{2 \pi s^3}} ds \nonumber \\
					& = \prod_{j=1}^{k-1} \lambda_j
					\sum_{m=1}^k
					\frac{1}{\prod_{
					l=1,l \neq m}^k \left( \lambda_l -
					\lambda_m \right) } \frac{\sin \nu \pi}{\pi} \int_0^\infty
					\frac{r^{\nu -1}}{r^{2 \nu} + 2 r^\nu \cos \nu \pi +1}
					e^{- t \sqrt{2 r \lambda_m^{\frac{1}{\nu}}}} dr \nonumber \\
					& = \frac{\sin \nu \pi}{\pi} \int_0^\infty
					\frac{r^{\nu -1}}{r^{2 \nu} + 2 r^\nu \cos \nu \pi +1}
					\prod_{j=1}^{k-1} \lambda_j
					\sum_{m=1}^k
					\frac{1}{\prod_{
					l=1,l \neq m}^k \left( \lambda_l -
					\lambda_m \right) }
					e^{- t \sqrt{2 r \lambda_m^{\frac{1}{\nu}}}} dr . \nonumber
				\end{align}
				If $\nu =1/2$ we obtain the following expression
				\begin{align}
					\label{interp}
					\hat{\mathpzc{p}}_k^{\frac{1}{2}} (t) & =
					\frac{1}{\pi} \int_0^\infty \frac{1}{\sqrt{r} \left( r +1 \right)}
					\prod_{j=1}^{k-1} (\sqrt{2r} \lambda_j)
					\sum_{m=1}^k
					\frac{1}{\prod_{
					l=1,l \neq m}^k \left( \sqrt{2r} \lambda_l -
					\sqrt{2r} \lambda_m \right) }
					e^{- t \sqrt{2 r} \lambda_m} dr \\
					& = \frac{ \sqrt{2} }{ \pi } \int_0^\infty \frac{1}{ \left( \frac{w^2}{2}
					+1 \right) } \text{Pr} \left\{ \mathpzc{N}_w (t) = k \mid \mathpzc{N}_w
					(0) =1 \right\} dw \nonumber
				\end{align}
				where $\mathpzc{N}_W (t)$, $t>0$, is a classical nonlinear birth process
				\eqref{nlinear} with random birth rates
				$(W \lambda_k)$, $k \geq 1$ where $W$ is a folded Cauchy r.v.
				with p.d.f.
				\begin{equation}
					f_W (w) = \frac{\sqrt{2}}{\pi \left( \frac{w^2}{2} +1 \right)},
					\qquad w \in \mathbb{R}^+ .
				\end{equation}
				It is possible to highlight a further interpretation by rewriting
				formula \eqref{interp} in the following way
				\begin{align}
					\hat{\mathpzc{p}}_k^{\frac{1}{2}} (t) & =
					\frac{\sqrt{2}}{\pi} \int_0^\infty \frac{1}{t \left( \frac{s^2}{2t^2}
					+1 \right)}
					\prod_{j=1}^{k-1} \lambda_j
					\sum_{m=1}^k
					\frac{1}{\prod_{
					l=1,l \neq m}^k \left( \lambda_l -
					\lambda_m \right) }
					e^{- s \lambda_m} ds \\
					& = \int_0^\infty \text{Pr} \left\{ \mathpzc{N} (s) = k \right\}
					\text{Pr} \left\{ \left| C(\sqrt{2} t ) \right| \in ds \right\}
					\nonumber
				\end{align}
				where $C(\sqrt{2} t)$, $t>0$ is a Cauchy process with rescaled time, possessing
				transition density
				\begin{equation}
					f_C (t,s) = \frac{1}{\pi} \frac{ \sqrt{2} t}{s^2+
					\left( \sqrt{2} t \right)^2},
					\qquad t>0, \: s \in \mathbb{R} .
				\end{equation}
				The process $\hat{\mathpzc{N}}^{1/2} (t) = \mathpzc{N}^{1/2}
				(T_t)$ can thus be written as $\hat{\mathpzc{N}}^{1/2} (t) =
				\mathpzc{N} ( \left| C(\sqrt{2} t ) \right| )$.
			\end{rem}

			\subsubsection{Iterated compositions}

				In the next theorem we present
				the explicit form of the state probabilities
				$\tilde{\mathpzc{p}}_k^\nu (t)$, $t>0$, $k \geq 1$, for the process
				\begin{equation}
					\tilde{\mathpzc{N}}^\nu (t) =
					\mathpzc{N}^\nu \biggl[ T^1_{T^2_{._{._{{T^n_t} }}}} \biggr],
					\qquad t>0, \: \nu \in (0,1].
				\end{equation}
				and in the following remark an interesting
				interpretation for that process when $\nu=1/2^n$,
				$n \in \mathbb{N}$, is given.

				\begin{thm}
					Let $\mathpzc{N}^\nu (t)$, $t>0$, $\nu \in (0,1]$,
					be a fractional nonlinear pure birth process and
					let $T^1_t$, $T^2_t$, $\dots T^n_t$, $n \in \mathbb{N}$,
					be $n$ independent first-passage time processes at $t$ of the
					standard Brownian motion.
					The process
					\begin{equation}
						\tilde{\mathpzc{N}}^\nu (t) =
						\mathpzc{N}^\nu \biggl[ T^1_{T^2_{._{._{{T^n_t} }}}} \biggr],
						\qquad t>0, \: \nu \in (0,1]
					\end{equation}
					has the following distribution
					\begin{equation}
						\label{nlinearnu-mcomp}
						\tilde{\mathpzc{p}}_k^\nu (t) =
						\begin{cases}
							\prod_{j=1}^{k-1} \lambda_j
							\sum_{m=1}^k \frac{1
							}{\prod_{
							l=1,l \neq m}^k \left( \lambda_l -
							\lambda_m \right) }
							\frac{\sin \nu \pi}{\pi}
							\int_0^\infty
							\frac{r^{\nu -1} e^{-t r^{\frac{1}{2^n}}
							\lambda_m^{\frac{1}{\nu 2^n}}
							2^{\left( 1- \frac{1}{2^n} \right)} } }{
							r^{2 \nu} + 2 r^\nu \cos \nu \pi +1}
							dr, & k > 1, \: t>0, \\
							\frac{\sin \nu \pi}{\pi}
							\int_0^\infty
							\frac{r^{\nu -1}
							e^{-t \lambda_1^{\frac{1}{
							2^n}} 2^{\left(1-\frac{1}{
							2^n} \right)}}}{
							r^{2 \nu} + 2 r^\nu \cos \nu \pi +1}
							dr, & k=1, \: t>0.
						\end{cases}
					\end{equation}
					\begin{proof}
						We start by proving the case $n=2$ since the
						case $n=1$ is already proved in Theorem \ref{teor}.
						We omit the details for the case $k=1$ and directly treat the
						case $k \geq 2$. We have
						\begin{align}
							\label{compformnu}
							& \int_0^\infty \hat{\mathpzc{p}}_k^\nu (s) q(t,s) ds \\
							& = \prod_{j=1}^{k-1} \lambda_j
							\sum_{m=1}^k \frac{1
							}{\prod_{
							l=1,l \neq m}^k \left( \lambda_l -
							\lambda_m \right) }
							\frac{\sin \nu \pi}{\pi} \int_0^\infty
							\frac{r^{\nu -1}}{r^{2 \nu} +2 r^\nu
							\cos \nu \pi +1} \int_0^\infty
							e^{-s \left( 2r \right)^{\frac{1}{2}} \lambda_m^{
							\frac{1}{2 \nu}}}
							\frac{t e^{- \frac{t^2}{2 s}}}{\sqrt{
							2 \pi s^3}} ds \: dr \nonumber \\
							& = \prod_{j=1}^{k-1} \lambda_j
							\sum_{m=1}^k
							\frac{1}{\prod_{
							l=1,l \neq m}^k \left( \lambda_l -
							\lambda_m \right)}
							\frac{\sin \nu \pi}{\pi} \int_0^\infty
							\frac{r^{\nu -1}}{r^{2 \nu} + 2 r^\nu
							\cos \nu \pi +1}
							e^{- t \sqrt{2 \sqrt{2 r \lambda_m^{\frac{1}{\nu}}}}}
							dr . \nonumber
						\end{align}
						It is now straightforward to generalise
						formula \eqref{compformnu}
						for $n$ compositions, as follows
						\begin{align}
							\label{more}
							\tilde{\mathpzc{p}}_k^\nu (t) & =
							\prod_{j=1}^{k-1} \lambda_j
							\sum_{m=1}^k
							\frac{1}{\prod_{
							l=1,l \neq m}^k \left( \lambda_l -
							\lambda_m \right)}
							\frac{\sin \nu \pi}{\pi} \int_0^\infty
							\frac{r^{\nu -1}}{r^{2 \nu} + 2 r^\nu
							\cos \nu \pi +1}
							e^{- t r^{\frac{1}{2^n}}
							\lambda_m^{\frac{1}{\nu 2^n}} 2^{\sum_{i=1}^n
							\frac{1}{2^i}}} dr\\
							& = \prod_{j=1}^{k-1} \lambda_j
							\sum_{m=1}^k
							\frac{1}{\prod_{
							l=1,l \neq m}^k \left( \lambda_l -
							\lambda_m \right)}
							\frac{\sin \nu \pi}{\pi} \int_0^\infty
							\frac{r^{\nu -1}}{r^{2 \nu} + 2 r^\nu
							\cos \nu \pi +1}
							e^{- t r^{\frac{1}{2^n}}
							\lambda_m^{\frac{1}{\nu 2^n}} 2^{\left(
							1- \frac{1}{2^n} \right)}} dr . \nonumber
						\end{align}
					\end{proof}
				\end{thm}
				\begin{rem}
					Analogously to Remark \ref{interpretation}, for
					$\nu = 1/2^n$, $n \in \mathbb{N}$, it is possible
					to interpret formula \eqref{nlinearnu-mcomp} as follows
					\begin{align}
						& \tilde{\mathpzc{p}}_k^{\frac{1}{2^n}} (t)
						= \prod_{j=1}^{k-1} \lambda_j
						\sum_{m=1}^k
						\frac{1}{\prod_{
						l=1,l \neq m}^k \left( \lambda_l -
						\lambda_m \right)}
						\frac{\sin \frac{\pi}{2^n}}{\pi} \int_0^\infty
						\frac{r^{\frac{1}{2^n} -1}}{r^{\frac{1}{2^{n -1}}} + 2
						r^{\frac{1}{2^n}} \cos \frac{\pi}{2^n} +1}
						e^{- t r^{\frac{1}{2^n}}
						2^{\left( 1- \frac{1}{2^n} \right)}
						\lambda_m } dr \\
						& \left(r=y^{2^n}\right) \notag \\
						& = \frac{\sin\frac{\pi}{2^n}}{\frac{\pi}{2^n}}
						\int_0^\infty \frac{dr}{r^2 +2r \cos \frac{\pi}{2^n}
						+1} \text{Pr} \left\{ \mathpzc{N}\left(
						tr2^{1-\frac{1}{2^n}} \right) = k \right\}. \notag
					\end{align}
					Therefore, the following representation holds
					\begin{equation}
						\mathpzc{N}^\frac{1}{2^n}
						\biggl[ T^1_{T^2_{._{._{{T^n_t} }}}} \biggr]
						= \mathpzc{N} \left( t\, \Omega\, 2^{1-\frac{1}{2^n}} \right),
					\end{equation}
					where $\Omega$ is a random variable with density
					\begin{equation}
						f_\Omega (r) =
						\frac{\sin \frac{\pi}{2^n}}{\frac{\pi}{2^n}}
						\frac{1}{r^2+2r\cos \frac{\pi}{2^n}+1},
						\qquad r \in \mathbb{R}^+ .
					\end{equation}
					The density is a unimodal law which, for $n\rightarrow\infty$,
					becomes
					\begin{equation}
						f(r) = \frac{1}{(1+r)^2}, \qquad r\in \mathbb{R}^+.
					\end{equation}
				\end{rem}

		\subsection{Fractional pure birth process stopped at $\mathpzc{S}^\alpha(t)$}

			We consider the fractional nonlinear pure birth process stopped at a
			stable time $\mathpzc{S}^\alpha(t)$ of order $0<\alpha\leq 1$ with Laplace
			transform
			\begin{equation}
				\mathbb{E} e^{-z \mathpzc{S}^\alpha(t)} = \int_0^\infty e^{-zs} q_\alpha
				(t,s) ds = e^{-tz^\alpha},
			\end{equation}
			where $q_\alpha(t,s)$, $s>0$, is the density of the stable process $\mathpzc{S}^\alpha(t)$,
			$t>0$.

			We have that the probabilities
			\begin{align}
				\label{compstanew}
				\breve{\mathpzc{p}}_k^\nu (t) & = \text{Pr} \left\{\mathpzc{N}^\nu
				(\mathpzc{S}^\alpha(t)) = k \mid \mathpzc{N}^\nu (0) = 1 \right\} \\
				& = \prod_{j=1}^{k-1} \lambda_j
				\sum_{m=1}^k
				\frac{1}{\prod_{
				l=1,l \neq m}^k \left( \lambda_l -
				\lambda_m \right) } \int_0^\infty E_{\nu,1} ( - \lambda_m s^\nu )
				q_\alpha(t,s) ds \nonumber \\
				& = \prod_{j=1}^{k-1} \lambda_j
				\sum_{m=1}^k
				\frac{1}{\prod_{
				l=1,l \neq m}^k \left( \lambda_l -
				\lambda_m \right) } \frac{\sin \nu \pi}{\pi} \int_0^\infty
				\frac{r^{\nu-1}}{r^{2 \nu} + 2r^\nu \cos \nu \pi +1}
				\int_0^\infty e^{-r \lambda_m^{\frac{1}{\nu}} s }
				q_\alpha(t,s) ds \, dr \nonumber \\
				& = \prod_{j=1}^{k-1} \lambda_j
				\sum_{m=1}^k
				\frac{1}{\prod_{
				l=1,l \neq m}^k \left( \lambda_l -
				\lambda_m \right) } \frac{\sin \nu \pi}{\pi} \int_0^\infty
				\frac{r^{\nu -1}}{r^{2 \nu} + 2 r^\nu \cos \nu \pi +1}
				e^{- t \lambda_m^{\frac{\alpha}{\nu}}r^\alpha} dr \nonumber \\
				& = \prod_{j=1}^{k-1} \lambda_j
				\sum_{m=1}^k
				\frac{1}{\prod_{
				l=1,l \neq m}^k \left( \lambda_l -
				\lambda_m \right) } \frac{\sin \nu \pi}{\alpha \pi}
				\int_0^\infty
				\frac{w^{\frac{\nu}{\alpha}-1}}{
				w^{2 \frac{\nu}{\alpha}} + 2 w^{\frac{\nu}{\alpha}} \cos \nu \pi
				+1} e^{- t \lambda_m^{\frac{\alpha}{\nu}} w} dw \nonumber \\
				& = \prod_{j=1}^{k-1} \lambda_j
				\sum_{m=1}^k \frac{1}{\prod_{l=1,l\neq m}^k \left( \lambda_l -\lambda_m \right)}
				\mathbb{E} e^{-t \lambda_m^{\frac{\alpha}{\nu}} \mathcal{W}_\alpha}, \notag
			\end{align}
			where $\mathcal{W}_\alpha$ is a random variable with density
			\begin{equation}
				\label{adens}
				f_{\mathcal{W}_\alpha} (w) = \frac{\sin \nu \pi}{\alpha \pi}
				\frac{w^{\frac{\nu}{\alpha}-1}}{
				w^{2 \frac{\nu}{\alpha}} + 2 w^{\frac{\nu}{\alpha}} \cos \nu \pi
				+1}, \qquad w>0, \: 0 < \nu < 1,
			\end{equation}
			(first obtained by \citet{lamperti}).
			The density \eqref{adens} coincides with the probability distribution of
			\begin{equation}
				\mathcal{W}_\alpha = \left( \frac{S_1^\nu}{S_2^\nu} \right)^\alpha,
			\end{equation}
			where $S_1^\nu$, $S_2^\nu$, are independent stable random variables with
			Laplace transform
			\begin{equation}
				\mathbb{E} e^{-z S^\nu} = e^{-z^\nu}, \qquad z>0, \: 0<\nu<1.
			\end{equation}

			If $\alpha=1$, $\mathbb{E} e^{-t \lambda_m^{\frac{1}{\nu}} \mathcal{W}_\alpha}
			= E_{\nu,1}(-\lambda_mt^\nu)$ and \eqref{compstanew} are the state
			probabilities of a fractional pure birth process, while for $\alpha=\nu$,
			$\mathbb{E} e^{-t \lambda_m \mathcal{W}_\nu}$ are the state probabilities
			of a pure birth process at time $t\mathcal{W}_\nu$ or, equivalently,
			a pure birth process at time $t$ with rates $\lambda_k \mathcal{W}_\nu$.

			If we compare \eqref{compstanew} with \eqref{nlinear-sub}, we can conclude that
			the process
			\begin{equation}
				\mathpzc{N}^\nu (\mathpzc{S}^\alpha(t)) =
				\mathpzc{N}(T_{2\nu}(\mathpzc{S}^\alpha(t))),
			\end{equation}
			can be represented as
			\begin{equation}
				\mathpzc{N} \left(\mathpzc{S}^{\alpha/\nu}
				\left({t\mathcal{W}_\alpha}\right)\right),
			\end{equation}
			if $0<\alpha<\nu<1$.

			\begin{rem}
				From formula \eqref{compstanew},
				when $\alpha$ takes the form $\alpha=\nu / 2^n$, $n \in \mathbb{N}$,
				we have
				\begin{equation}
					\breve{\mathpzc{p}}_k^\nu (t) = \frac{\sin \nu \pi}{\pi} \int_0^\infty
					\frac{r^{\nu -1}}{r^{2 \nu} + 2 r^\nu \cos \nu \pi +1}
					\prod_{j=1}^{k-1} \lambda_j
					\sum_{m=1}^k
					\frac{e^{-t \lambda_m^{\frac{1}{2^n}}
					r^{\frac{\nu}{2^n}}}}{\prod_{
					l=1,l \neq m}^k \left( \lambda_l -
					\lambda_m \right) } dr .
				\end{equation}
				For $n \rightarrow \infty$, we obtain that $\breve{\mathpzc{p}}_k^\nu (t)
				\rightarrow 0$, $k > 1$, and $\breve{\mathpzc{p}}_1^\nu (t) \rightarrow
				e^{-t}$.
				This shows that for $n\rightarrow \infty$, either the population
				istantaneously explodes or does not produce offsprings with
				exponential probability.
			\end{rem}

			An alternative way of presenting the state probabilities \eqref{compstanew}
			is based on the Mellin--Barnes representation of the Mittag--Leffler
			function
			\begin{equation}
				\label{barnes}
				E_{\nu,\mu} (x) = \frac{1}{2 \pi i}\int_{\gamma-i\infty}^{\gamma+i\infty}
				\frac{\Gamma(z) \Gamma(1-z)}{\Gamma(\mu- \nu z)} x^{-z} dz,
			\end{equation}
			with $\nu>0$, $x \in \mathbb{C}$, $\left| \arg (-x) \right| < \pi$
			(see \citet{kilbas}, page 44, formula (1.8.32)).

			In view of \eqref{barnes}, we can write \eqref{compstanew} as follows
			\begin{align}
				\label{rewri}
				\breve{\mathpzc{p}}_k^\nu (t) & =
				\prod_{j=1}^{k-1} \lambda_j
				\sum_{m=1}^k
				\frac{1}{\prod_{
				l=1,l \neq m}^k \left( \lambda_l -
				\lambda_m \right) } \int_0^\infty
				\frac{1}{2 \pi i} \int_{\gamma -i \infty}^{\gamma +i \infty}
				\frac{\Gamma(z) \Gamma(1-z)}{\Gamma(1- \nu z)}
				(\lambda_m s^\nu)^{-z} dz \, q_\alpha(t,s) ds \\
				& = \prod_{j=1}^{k-1} \lambda_j
				\sum_{m=1}^k
				\frac{1}{\prod_{
				l=1,l \neq m}^k \left( \lambda_l -
				\lambda_m \right) } \frac{1}{2 \pi i}
				\int_{\gamma -i \infty}^{\gamma +i \infty}
				\frac{\Gamma(z) \Gamma(1-z)}{\Gamma(1- \nu z)}
				\lambda_m^{-z} \int_0^\infty s^{-\nu z} q_\alpha (t,s) ds \, dz,
				\notag
			\end{align}
			where in the last member of \eqref{rewri} the Mellin transform of
			$q_\alpha(t,s)$ appears.

			The Mellin transform of the stable subordinator $\mathpzc{S}^\alpha(t)$, with
			Laplace transform
			\begin{equation}
				\mathbb{E}e^{-z\mathpzc{S}^\alpha(t)} = e^{-tz^\alpha},
			\end{equation}
			reads
			\begin{align}
				\label{nillem}
				\mathbb{E}(\mathpzc{S}^\alpha(t))^{\eta-1} & = \int_0^\infty
				s^{\eta-1} q_\alpha(t,s) ds
				= \frac{1}{\alpha} \Gamma \left( \frac{1-\eta}{\alpha} \right)
				\frac{1}{\Gamma(1-\eta)} t^{\frac{\eta-1}{\alpha}}.
			\end{align}

			By inserting \eqref{nillem} into \eqref{rewri}, we arrive at
			\begin{align}
				\label{llll}
				\breve{\mathpzc{p}}_k^\nu (t) & =
				\prod_{j=1}^{k-1} \lambda_j
				\sum_{m=1}^k
				\frac{1}{\prod_{
				l=1,l \neq m}^k \left( \lambda_l -
				\lambda_m \right) } \frac{\alpha^{-1}}{2 \pi i}
				\int_{\gamma -i\infty}^{\gamma + i \infty}
				\frac{ \Gamma(z) \Gamma(1-z) \Gamma
				\left( \frac{\nu}{\alpha}z \right)}{
				\Gamma( \nu z) \Gamma (1-\nu z)} \left( \lambda_m t^{
				\frac{\nu}{\alpha}} \right)^{-z} dz \\
				& = \prod_{j=1}^{k-1} \lambda_j
				\sum_{m=1}^k
				\frac{1}{\prod_{
				l=1,l \neq m}^k \left( \lambda_l -
				\lambda_m \right) } \alpha^{-1}
				H^{2,1}_{2,3} \left[ \lambda_m t^{\frac{\nu}{\alpha}}
				\left|
				\begin{array}{c}
					(0,1), (0,\nu ) \\
					(0,1), (0,\nu/\alpha), (0,\nu)
				\end{array}
				\right. \right]. \notag
			\end{align}

			We examine now in detail the case $\alpha=\nu$ in the next theorem.
			\begin{thm}
				We have the following distributions:
				\begin{enumerate}
					\item $\text{Pr} \left\{ \mathpzc{N}^\nu(\mathpzc{S}^\nu(t)) = k \right\}
						= \text{Pr} \left\{ \mathpzc{N} (T_{2\nu}(\mathpzc{S}^\nu(t)))
						= k \right\} = \text{Pr} \left\{ \mathpzc{N} (t \mathcal{W}_\nu)
						= k \right\}$,
					\item $\text{Pr} \left\{ \mathpzc{N}
						\left(\mathpzc{S}^\nu({T_{2\nu}(t)}) \right) = k
						\right\} = \text{Pr} \left\{ \mathpzc{N} (t \mathcal{W}_1) = k
						\right\}$,
				\end{enumerate}
				for $k\geq 1$, $t>0$, where
				\begin{equation}
					\mathcal{W}_\alpha = \left( \frac{S_1^\nu}{S_2^\nu} \right)^\alpha,
				\end{equation}
				and has distribution \eqref{adens}.

				\begin{proof}
					For $k>1$ we can write that
					\begin{align}
						& \text{Pr} \left\{ \mathpzc{N} (T_{2\nu} ( \mathpzc{S}^\nu(t)
						)) = k \right\} \\
						& = \prod_{j=1}^{k-1} \lambda_j
						\sum_{m=1}^k
						\frac{1}{\prod_{
						l=1,l \neq m}^k \left( \lambda_l -
						\lambda_m \right) } \int_0^\infty
						e^{-\lambda_m s} \text{Pr} \left\{ T_{2\nu} (
						\mathpzc{S}^\nu(t)) \in ds \right\} \notag \\
						& = \prod_{j=1}^{k-1} \lambda_j
						\sum_{m=1}^k
						\frac{1}{\prod_{
						l=1,l \neq m}^k \left( \lambda_l -
						\lambda_m \right) } \int_0^\infty
						e^{-\lambda_m s} \int_0^\infty
						f_{T_{2\nu}} (z,s) q_\nu (t,z) dz \, ds \notag \\
						& = \prod_{j=1}^{k-1} \lambda_j
						\sum_{m=1}^k
						\frac{1}{\prod_{
						l=1,l \neq m}^k \left( \lambda_l -
						\lambda_m \right) } \int_0^\infty
						E_{\nu,1}(-\lambda_m z^\nu) q_\nu(t,z) dz \notag \\
						& = \prod_{j=1}^{k-1} \lambda_j
						\sum_{m=1}^k
						\frac{1}{\prod_{
						l=1,l \neq m}^k \left( \lambda_l -
						\lambda_m \right) } \frac{\sin \nu \pi}{\pi}
						\int_0^\infty dr \int_0^\infty
						\frac{r^{\nu-1}e^{-\lambda_m^{1/\nu} zr }}{
						r^{2\nu} +2r^\nu \cos \nu \pi +1 } q_\nu (t,z) dz \notag \\
						& = \prod_{j=1}^{k-1} \lambda_j
						\sum_{m=1}^k
						\frac{1}{\prod_{
						l=1,l \neq m}^k \left( \lambda_l -
						\lambda_m \right) } \frac{\sin \nu \pi}{\pi}
						\int_0^\infty \frac{r^{\nu-1} e^{-t(r\lambda_m^{1/\nu} )^\nu } }{
						r^{2\nu} +2r^\nu \cos \nu \pi +1 } dr \notag \\
						& = \int_0^\infty \text{Pr} \left\{ \mathpzc{N} (tr^\nu)=k
						\right\} \frac{\sin \nu \pi}{\pi} \frac{r^{\nu-1}}{
						r^{2\nu} +2r^\nu \cos \nu \pi +1 } dr \notag \\
						& = \text{Pr} \left\{ \mathpzc{N} ( t\mathcal{W}_1^\nu ) = k
						\right\} = \text{Pr} \left\{ \mathpzc{N} ( t\mathcal{W}_\nu)
						= k \right\}. \notag
					\end{align}
					This concludes the proof of the first result. In order to prove
					the second result we write
					\begin{align}
						& \text{Pr} \left\{ \mathpzc{N} \left( \mathpzc{S}^\nu({
						T_{2\nu}(t)}) \right) = k \right\} \\
						& = \prod_{j=1}^{k-1} \lambda_j
						\sum_{m=1}^k
						\frac{1}{\prod_{
						l=1,l \neq m}^k \left( \lambda_l -
						\lambda_m \right) } \int_0^\infty \int_0^\infty
						e^{-\lambda_m s} q_\nu (z,s) f_{T_{2\nu}} (z,t) dz \, ds
						\notag \\
						& = \prod_{j=1}^{k-1} \lambda_j
						\sum_{m=1}^k
						\frac{1}{\prod_{
						l=1,l \neq m}^k \left( \lambda_l -
						\lambda_m \right) } \int_0^\infty e^{-\lambda_m^\nu z}
						f_{T_{2\nu}} (z,t) dz \notag \\
						& \overset{\text{by} \eqref{ackard}}{=}
						\prod_{j=1}^{k-1} \lambda_j
						\sum_{m=1}^k
						\frac{1}{\prod_{
						l=1,l \neq m}^k \left( \lambda_l -
						\lambda_m \right) } E_{\nu,1} (-\lambda_m^\nu t^\nu ) \notag \\
						& = \int_0^\infty
						\prod_{j=1}^{k-1} \lambda_j
						\sum_{m=1}^k
						\frac{1}{\prod_{
						l=1,l \neq m}^k \left( \lambda_l -
						\lambda_m \right) } \frac{\sin \nu \pi}{\pi}
						\frac{r^{\nu-1} e^{-\lambda_m t r}}{r^{2\nu} +
						2r^\nu \cos \nu \pi + 1} dr \notag \\
						& = \int_0^\infty \text{Pr} \left\{ \mathpzc{N}
						(tr) = k \right\} \frac{\sin \nu \pi}{\pi} \frac{r^{\nu-1}}{
						r^{2\nu} +2r^\nu \cos \nu \pi +1 } dr
						= \text{Pr} \left\{ \mathpzc{N} (t \mathcal{W}_1) = k \right\}
						. \notag
					\end{align}
				\end{proof}
			\end{thm}

			\begin{rem}
				By slightly changing the above calculations, we arrive at the
				following result (compare with \eqref{compstanew}):
				\begin{equation}
					\text{Pr} \left\{ \mathpzc{N} \left( \mathpzc{S}^\alpha({
						T_{2\nu}(t)}) \right) = k \right\}
						= \prod_{j=1}^{k-1} \lambda_j
						\sum_{m=1}^k
						\frac{1}{\prod_{
						l=1,l \neq m}^k \left( \lambda_l -
						\lambda_m \right) } \mathbb{E}
						e^{-t\lambda_m^{\frac{\alpha}{\nu}} \mathcal{W}_1 }.
				\end{equation}
			\end{rem}

			\begin{rem}
				An alternative form of the distribution \eqref{compstanew}, for $\alpha=\nu$,
				can be given as follows.
				\begin{align}
					& \text{Pr} \left\{ \mathpzc{N}^\nu ( \mathpzc{S}^\nu(t)) = k \right\} \\
					& = \prod_{j=1}^{k-1} \lambda_j
					\sum_{m=1}^k
					\frac{1}{\prod_{
					l=1,l \neq m}^k \left( \lambda_l -
					\lambda_m \right) } \frac{1}{2 i \nu \pi} \int_0^\infty
					\left( \frac{1}{w+e^{-i \pi \nu}} - \frac{1}{w+e^{i \pi \nu}}
					\right) e^{-t \lambda_m w} dw \notag \\
					& = \prod_{j=1}^{k-1} \lambda_j
					\sum_{m=1}^k
					\frac{1}{\prod_{
					l=1,l \neq m}^k \left( \lambda_l -
					\lambda_m \right) } \frac{1}{2 i \nu \pi}
					\left[ e^{t \lambda_m e^{-i \pi \nu}} \mathrm{E}_1 \left( t \lambda_m
					e^{-i \pi \nu} \right) - e^{t \lambda_m e^{i \pi \nu}} \mathrm{E}_1
					\left( t \lambda_m e^{i \pi \nu} \right) \right] \notag \\
					& = \prod_{j=1}^{k-1} \lambda_j
					\sum_{m=1}^k
					\frac{1}{\prod_{
					l=1,l \neq m}^k \left( \lambda_l -
					\lambda_m \right) } \frac{ e^{\lambda_m t \cos \nu \pi} }{2 i \nu \pi}
					\left[ e^{t \lambda_m \sin \nu \pi} \mathrm{E}_1 \left( t \lambda_m
					e^{-i \pi \nu} \right) - e^{t \lambda_m \sin \nu \pi} \mathrm{E}_1
					\left( t \lambda_m e^{i \pi \nu} \right) \right], \notag
				\end{align}
				where the function
				$\mathrm{E}_1 (z) = \int_z^\infty \frac{e^{-t}}{t} dt$,
				$|\arg z< \pi|$, is the exponential integral.
			\end{rem}

		\subsection{Fractional pure birth process stopped at $T_{2\alpha}(t)$}

			In this section we consider the process $\mathpzc{N}^\nu( T_{2 \alpha} (t))$, $t>0$
			(see the discussion related to formula \eqref{diffusion}. As we did before,
			here we treat the case $k \geq 2$.
			The state probabilities can be written as follows.
			\begin{align}
				\label{ricomp}
				\mathpzc{p}_k^{\nu,\alpha} (t) & = \text{Pr} \left\{\mathpzc{N}^\nu
				(T_{2\alpha}(t)) = k \mid \mathpzc{N}^\nu (0) = 1 \right\} \\
				& = \prod_{j=1}^{k-1} \lambda_j
				\sum_{m=1}^k
				\frac{1}{\prod_{
				l=1,l \neq m}^k \left( \lambda_l -
				\lambda_m \right) } \int_0^\infty E_{\nu,1} ( - \lambda_m s^\nu )
				\text{Pr} \left\{T_{2\alpha} (t) \in ds \right\}. \notag
			\end{align}
			The integral in \eqref{ricomp} can be further worked out by means of the Laplace
			transform:
			\begin{align}
				& \int_0^\infty E_{\nu,1} (-\lambda_m s^\nu) \int_0^\infty e^{-zt}
				\text{Pr} \left\{T_{2\alpha} (t) \in ds \right\} dt \\
				& = \int_0^\infty E_{\nu,1} (-\lambda_m s^\nu) z^{\alpha-1}
				e^{-sz^\alpha} ds
				= z^{\alpha -1} \frac{(z^\alpha)^{\nu-1}}{z^{\alpha \nu}+\lambda_m}
				= \frac{z^{\alpha \nu -1}}{z^{\alpha \nu} + \lambda_m}. \notag
			\end{align}
			By taking the inverse Laplace transform of the above formula, we immediately obtain that
			\begin{align}
				\int_0^\infty E_{\nu,1} (-\lambda_m s^\nu) \text{Pr} \left\{
				T_{2\alpha} (t) \in ds \right\}
				= \int_0^\infty E_{\nu,1} (-\lambda_m s^\nu)
				t^{-\alpha} W_{-\alpha,1-\alpha} (-t^{-\alpha} s) ds
				= E_{\nu \alpha, 1} (-\lambda_m t^{\nu \alpha}).
			\end{align}
			Therefore, the state probabilities for the process $\mathpzc{N}^\nu(
			T_{2 \alpha} (t))$, $t>0$, result in the following form:
			\begin{equation}
				\label{expl2}
				\mathpzc{p}_k^{\nu,\alpha} (t) =
				\prod_{j=1}^{k-1} \lambda_j
				\sum_{m=1}^k
				\frac{1}{\prod_{
				l=1,l \neq m}^k \left( \lambda_l -
				\lambda_m \right) } E_{\nu \alpha, 1} (-\lambda_m t^{\nu \alpha})
				= \mathpzc{p}_k^{\nu \alpha}(t), \qquad k \geq 2.
			\end{equation}
			Note that the case $k=1$ can be treated in the same manner.
			We thus obtain the following equalities in distribution:
			\begin{equation}
				\label{formula2}
				\mathpzc{N}^\nu(T_{2\alpha}(t)) = \mathpzc{N}\left\{T_{2\nu}
				(T_{2\alpha}(t))\right\} = \mathpzc{N} (T_{2\nu\alpha}(t))
				= \mathpzc{N}^{\nu \alpha}(t), \qquad t>0.
			\end{equation}
			Let now $\eta_n = \prod_{i=1}^n \nu_i$, where $n \in \mathbb{N}$, and $\nu_i$
			are $n$ indices such that $\nu_i \in (0,1]$ for $1\leq i\leq n$.
			Formula \eqref{formula2} can be generalised as
			\begin{equation}
				\mathpzc{N} \left\{ T_{2\nu_1} ( T_{2\nu_2} (\dots T_{2\nu_n}(t)\dots ))
				\right\} = \mathpzc{N} (T_{2\eta_n}(t)) = \mathpzc{N}^{\eta_n} (t),
				\qquad t>0
			\end{equation}
			where $\mathpzc{N}^{\eta_n}(t)$ is a nonlinear fractional birth process.

			Formula \eqref{ricomp} can also be worked out in an alternative way.
			In the following calculations we will make use of the integral
			representation \eqref{mitta-int}.
			\begin{align}
				\label{ricomp2}
				\mathpzc{p}_k^{\nu,\alpha} (t) = {} & \text{Pr} \left\{\mathpzc{N}^\nu
				(T_{2\alpha}(t)) = k \mid \mathpzc{N}^\nu (0) = 1 \right\} \\
				= {} & \prod_{j=1}^{k-1} \lambda_j
				\sum_{m=1}^k
				\frac{1}{\prod_{
				l=1,l \neq m}^k \left( \lambda_l -
				\lambda_m \right) } \int_0^\infty E_{\nu,1} \left( - \lambda_m s^\nu \right)
				\text{Pr} \left\{T_{2\alpha} (t) \in ds \right\} \notag \\
				= {} & \int_0^\infty \prod_{j=1}^{k-1} \lambda_j
				\sum_{m=1}^k
				\frac{1}{\prod_{
				l=1,l \neq m}^k \left( \lambda_l -
				\lambda_m \right) } \notag \\
				& \times \int_0^\infty \frac{\sin \nu \pi}{\pi}
				\frac{r^{\nu-1}}{r^{2\nu}+2r^\nu \cos \nu \pi+1}
				e^{-r \lambda_m^{\frac{1}{\nu}} s}
				\text{Pr} \left\{T_{2\alpha} (t) \in ds \right\} dr \notag \\
				= {} & \prod_{j=1}^{k-1} \lambda_j
				\sum_{m=1}^k
				\frac{1}{\prod_{
				l=1,l \neq m}^k \left( \lambda_l -
				\lambda_m \right) } \int_0^\infty \frac{\sin \nu \pi}{\pi}
				\frac{r^{\nu-1}}{r^{2\nu}+2r^\nu \cos \nu \pi+1}
				E_{\alpha,1} \left(-r \lambda_m^{\frac{1}{\nu}} t^\alpha\right). \notag
			\end{align}

			\begin{rem}
				By comparing formulae \eqref{expl2} and \eqref{ricomp2}, it is
				clear that the following expansion holds:
				\begin{equation}
					\label{gen-mitta}
					E_{\nu \alpha, 1} (-\lambda_m t^{\nu \alpha}) =
					\frac{\sin \nu \pi}{\pi} \int_0^\infty
					\frac{r^{\nu-1}}{r^{2\nu}+2r^\nu \cos \nu \pi+1}
					E_{\alpha,1} \left(-r \lambda_m^{\frac{1}{\nu}} t^\alpha\right) dr,
					\qquad \nu \in (0,1], \: \alpha \in (0,1].
				\end{equation}

				We give a direct proof of \eqref{gen-mitta} by applying the Laplace transform
				to both members.
				Of course
				\begin{equation}
					\int_0^\infty e^{-\mu t} E_{\nu \alpha,1} (-\lambda_m t^{\nu\alpha})
					dt = \frac{\mu^{\nu\alpha -1}}{\mu^{\nu\alpha}+\lambda_m}.
				\end{equation}
				Then, we must calculate the twofold integral
				\begin{align}
					& \frac{\sin \nu\pi}{\pi} \int_0^\infty e^{-\mu t} \int_0^\infty
					\frac{r^{\nu-1}}{r^{2\nu}+2r^\nu\cos\nu\pi+1} E_{\alpha,1}
					\left( -r \lambda_m^{\frac{1}{\nu}} t^\alpha \right) dt \, dr \\
					& = \frac{\sin \nu\pi}{\pi} \int_0^\infty
					\frac{r^{\nu-1}}{r^{2\nu}+2r^\nu\cos\nu\pi+1} \frac{\mu^{\alpha -1}}{
					\mu^\alpha +r\lambda_m^{\frac{1}{\nu}}} dr \notag \\
					& = \frac{\sin \nu\pi}{\pi} \int_0^\infty \int_0^\infty
					\frac{\mu^{\alpha-1}r^{\nu-1}e^{-w \left(\mu^\alpha +r\lambda_m^{\frac
					1 \nu}\right)}}{r^{2\nu}+2r^\nu\cos\nu\pi+1} dr \, dw \notag \\
					& \overset{\text{by \eqref{mitta-int}}}{=}
					\int_0^\infty \mu^{\alpha -1} e^{-\mu^\alpha w}
					E_{\nu,1} (-\lambda_m w^\nu) dw
					= \frac{(\mu^\alpha)^{\nu-1} \mu^{\alpha-1}}{(\mu^\alpha)^\nu
					+\lambda_m} = \frac{\mu^{\nu\alpha-1}}{\mu^{\nu\alpha}+\lambda_m}. \notag
				\end{align}

				\begin{rem}
					A number of interesting relations follow from formula \eqref{gen-mitta}.

					The following integral relation holds:
					\begin{align}
						& \frac{\sin\nu\pi}{\nu} \int_0^\infty
						\frac{1}{r^2 +2r\cos\nu\pi +1} E_{\alpha,1}
						\left( -r^{\frac{1}{\nu}} \lambda_m^{\frac{1}{\nu}}
						t^\alpha \right) dr \\
						& = \frac{\sin\nu\alpha\pi}{\nu\alpha}
						\int_0^\infty \frac{1}{r^2 +2r\cos\nu\alpha\pi +1}
						e^{-r^{\frac{1}{\nu\alpha}} \lambda_m^{\frac{1}{\nu\alpha}}
						t} dr. \notag
					\end{align}

					A sort of commutativity is valid for \eqref{gen-mitta}:
					\begin{equation}
						\label{gen-mitta2}
						E_{\alpha \nu, 1} (-\lambda_m t^{\alpha \nu}) =
						\frac{\sin \alpha \pi}{\pi} \int_0^\infty
						\frac{r^{\alpha-1}}{r^{2\alpha}+2r^\alpha \cos \alpha \pi+1}
						E_{\nu,1} \left(-r \lambda_m^{\frac{1}{\alpha}} t^\nu\right),
						\qquad \nu \in (0,1], \: \alpha \in (0,1].
					\end{equation}

					For $\alpha=1$ we recover, from \eqref{gen-mitta}, the integral
					representation of Mittag-Leffler functions. By considering that
					\begin{equation}
						f_\nu(r) = \frac{r^{\nu-1}}{r^{2\nu}+2r^\nu \cos \nu \pi+1},
					\end{equation}
					for $\nu=1$, becomes a delta function with pole at $r=1$, we extract,
					from \eqref{gen-mitta}, an identity.
				\end{rem}

				Furthermore, it is worth noticing that formulae similar to \eqref{gen-mitta} can be
				derived by repeated applications of the same formula.
				For example we have:
				\begin{align}
					& E_{\nu \alpha \beta,1} \left( -\lambda_m t^{\nu \alpha \beta} \right)
					= \frac{\sin \nu \pi}{\pi} \int_0^\infty
					\frac{r^{\nu-1}}{r^{2\nu}+2r^\nu \cos \nu \pi+1} E_{\alpha \beta,1}
					\left( -r \lambda_m^{\frac{1}{\nu}} t^{\alpha \beta} \right) dr \\
					& = \frac{\sin \nu \pi \sin \alpha \pi}{\pi^2} \int_0^\infty
					\int_0^\infty \frac{r^{\nu-1}w^{\alpha-1}
					E_{\beta,1}( -w r^\frac{1}{\alpha} \lambda_m^{\frac{1}{\nu
					\alpha}} t^\beta )}{(r^{2\nu}+2r^\nu \cos
					\nu \pi+1)(w^{2\alpha}+2w^\alpha \cos \alpha \pi+1)}\, dw \, dr. \notag
				\end{align}
				Let $\nu_i$, $1\leq i\leq n$ be $n$ indices such that
				for all $1\leq i \leq n$, $\nu_i \in (0,1]$, and
				let us denote $\eta_n=\prod_{i=1}^n \nu_i$. In general, for $n \geq 2$, we obtain
				that
				\begin{align}
					& E_{\eta_n,1} \left( -\lambda_m t^{\eta_n}\right) \\
					& = \int_0^\infty \hspace{-0.35cm} \dots \int_0^\infty
					\prod_{j=1}^{n-1} \left( \frac{ r_j^{\nu_j-1} }{
					r_j^{2\nu_j}+2 r_j^{\nu_j}\cos \nu_j +1  } \right)
					E_{\nu_n,1} \left(
					- r_1 r_2^{\frac{1}{\nu_1}} r_3^{\frac{1}{\nu_1\nu_2}} \dots
					r_{n-1}^{\frac{1}{\eta_{n-2}}} \lambda_m^{\frac{1}{\eta_{n-1}}} t^{\nu_n}
					\right) \prod_{j=1}^{n-1} dr_j \notag \\
					& =\mathbb{E} \left[ E_{\nu_n,1} \left( - {}^{(1)}\mathcal{W}_1
					{}^{(2)}\mathcal{W}_1^{\frac{1}{\nu_1}}
					{}^{(3)}\mathcal{W}_1^{\frac{1}{\nu_1\nu_2}}
					\dots {}^{(n-1)}\mathcal{W}_1^{
					\frac{1}{\eta_{n-2}}} \lambda_m^{\frac{1}{\eta_{n-1}}} t^{\nu_n}
					\right) \right], \notag
				\end{align}
				where ${}^{(j)}\mathcal{W}_1$, $1\leq j \leq n-1$,
				are independent random variables,
				each with distribution \eqref{adens}, with $\alpha=1$ and $\nu=\nu_j$.
			\end{rem}

	\bibliographystyle{abbrvnat}
	\bibliography{composition3}
	\nocite{*}

\end{document}